\documentclass[reqno,centertags,a4paper]{amsart}
\usepackage{amssymb,amsmath,amsthm}
\usepackage{hyperref}
\usepackage{times}

\newcommand{\R}{\mathbb{R}}
\newcommand{\Z}{\mathbb{Z}}
\newcommand{\T}{\mathbb{T}}
\newcommand{\C}{\mathbb{C}}
\newcommand{\N}{\mathbb{N}}

\theoremstyle{plain}
\newtheorem{theorem}{Theorem}[section]
\newtheorem{lemma}{Lemma}[section]
\newtheorem{corollary}{Corollary}[section]

\theoremstyle{definition}
\newtheorem{definition}{Definition}[section]
\theoremstyle{remark}
\newtheorem{remark}{Remark}

\DeclareMathOperator{\supp}{supp}
\DeclareMathOperator{\sign}{sign}

\DeclareMathOperator{\Imag}{Im}
\DeclareMathOperator{\Real}{Re}

\newcommand{\diff}[1]{\frac{{\mathrm d}}{{\mathrm d #1}}}
\newcommand{\test}{\mathcal{S}_{\mathrm per}}
\newcommand{\eps}{\varepsilon}
\newcommand{\lb}{\langle}
\newcommand{\rb}{\rangle}

\author[S. Herr]{Sebastian Herr}

\title[On the Cauchy Problem for the DNLS with periodic boundary condition]
{On the Cauchy problem for the derivative nonlinear Schr\"odinger
equation with periodic boundary condition}


\address{Fachbereich Mathematik, Universit\"at Dortmund, 44221 Dortmund, Germany.}
\email{sebastian.herr@math.uni-dortmund.de}

\begin{document}

\begin{abstract}
It is shown that the Cauchy problem associated
to the derivative nonlinear Schr\"odinger equation
$\partial_t u-i\partial^2_{x}u=\lambda \partial_{x}(|u|^2u)$
is locally well-posed for initial data $u(0)\in
H^{s}(\T)$, if $s\geq \frac{1}{2}$ and $\lambda$ is real.
The proof is based on an adaption of the
gauge transformation to periodic functions and sharp multi-linear estimates
for the gauge equivalent equation in Fourier restriction norm spaces. By the use
of a
conservation law, the problem is shown to be globally well-posed for $s\geq 1$
and data which is small in $L^2$.
\end{abstract}


\maketitle

\section{Introduction and main result}\label{sect:main}
\noindent
We study the Cauchy problem associated to the derivative
nonlinear Schr\"odinger (DNLS) equation with the periodic boundary condition
\begin{equation}\label{eq:pdnls}
\begin{split}
\partial_t u-i\partial^2_{x}u&=\lambda \partial_{x}(|u|^2u) \quad\text{in } (-T,T)\times \T\\
u(0)&=u_0 \in H^{s}(\T)
\end{split}
\end{equation}
where $\T=\R/2\pi\Z$ and $\lambda \in \R$.
Our aim is to prove local and global well-posedness in low regularity Sobolev spaces.

In the case of the real line, local well-posedness in $H^s(\R)$ for $s\geq
\frac{1}{2}$ was obtained by Takaoka
\cite{Tak99} and this was shown to be sharp in the sense of the uniform continuity of
the flow map by Biagoni and Linares \cite{BL01} (the critical regularity for the
scaling argument is $L^2$). The local result was extended to global well-posedness
for $s >\frac{1}{2}$ by Colliander, Keel, Staffilani, Takaoka and Tao
\cite{CKSTT02} for data which satisfies a $L^2$ smallness condition.
Recently, Gr\"unrock \cite{G05} obtained a local result in a $\widehat{L^p}$
setting. The DNLS equation found application as a model in plasma physics
and it satisfies infinitely many conservation laws \cite{KN78}.
For a more detailed history and further references we refer the reader to these works.

The local result of Takaoka was proved by using the gauge transform developed by
Hayashi and Ozawa \cite{HO92,Hay93,HO94} to derive a gauge equivalent equation. This
equation still contains a tri-linear term with derivative of the form $u^2
\partial_x\overline{u}$, but Takaoka \cite{Tak99} showed that this can
be treated by the Fourier restriction norm method developed in \cite{Bo93}, as long as
$s\geq \frac{1}{2}$. The proof of the main tri-linear estimate
uses local smoothing and Strichartz estimates.

Here, we study the DNLS equation in the periodic setting.
It is known that there exist global (weak) solutions in Sobolev spaces corresponding
to $H^1$ subject to Dirichlet and generalized periodic boundary conditions due to the results from
Chen \cite{Chen86} and Me\v{s}kauskas \cite{Me98}, for initial data fulfilling a smallness condition.

We remark that the dispersive properties of solutions are weaker than in the
non-periodic case. Of course, there are no local smoothing estimates available which could
be used to control derivatives in nonlinear terms. Moreover, above $L^4$ the
Strichartz estimates are only known to hold with a loss of $\eps>0$ derivatives.
Therefore, the main question arising in the periodic case is whether a tri-linear estimate for $u^2
\partial_x\overline{u}$ holds true.

In the present work we will answer this question affirmatively. Our main
ingredients are a point-wise estimate for the multiplier,
suitable versions of Bourgain spaces \cite{Bo93,G96,GTV97,CKSTT03,G00} and the
$L^4$ Strichartz estimate \cite{Bo93}.
Combining this with a gauge transform \cite{HO92,Hay93,HO94} adapted to the
periodic setting, it follows that
the Cauchy problem \eqref{eq:pdnls} is locally well-posed in $H^s(\T)$ for $s\geq
\frac{1}{2}$ in the following sense:
\begin{theorem}\label{thm:main}
Let $s\geq \frac{1}{2}$ and $\lambda \in \R$. For all $r>0$ there exists
$T=T(r)>0$ and a metric space
$\mathcal{X}_s^T$, such that for all $ u_0 \in B_{r}=\{u_0 \in H^{s}(\T)\mid
\|u_0\|_{H^\frac{1}{2}(\T)}< r\}$
there exists a unique solution
$$u \in \mathcal{X}_s^T \hookrightarrow C([-T,T],H^s(\T))$$
of the Cauchy problem
\begin{align*}
\partial_t u-i\partial^2_{x}u&=\lambda \partial_{x}(|u|^2u) \quad\text{in } (-T,T)\times \T\\
u(0)&=u_0
\end{align*}
which is a limit of smooth solutions in $\mathcal{X}_s^T$.

Moreover, the flow map 
$$
F: H^s(\T) \supset B_{r}
\to C\big([-T,T],H^{s}(\T)\big) \quad,\; u_0 \mapsto u
$$
is continuous. For fixed $\mu > 0$ the restriction of $F$ to $\{u_0
\in B_r \mid \frac{1}{2\pi}\|u_0\|^2_{L^2}=\mu\}$
is Lipschitz continuous.
\end{theorem}
Due to a conservation law this extends to global
well-posedness for $s\geq 1$ for data which is small in $L^2$.
\begin{corollary}\label{cor:global}
Let $s\geq 1$. There exists $\delta>0$ such that under the additional hypothesis
$\|u_0\|_{L^2}\leq\delta$, the time of existence $T>0$ in Theorem \ref{thm:main} can
be chosen arbitrary large.
\end{corollary}
Throughout this work \emph{solution} of a Cauchy problem
\begin{align*}
\partial_t u-i\partial_x^2 u &=N(u) \quad\text{in } (-T,T)\times \T\\
u(0)&=u_0
\end{align*}
always means solution of the corresponding integral equation
\begin{equation*}
u(t)=W(t)u_0+\int_0^t W(t-t') N(u)(t')\,dt', \quad t \in (-T,T)
\end{equation*}
at least in a limiting sense,
see Sections \ref{sect:thm_gauge_eq_problem} and \ref{sect:proof_thm_main}
for the precise statements.
For smooth functions, this notion of solutions coincides with the classical one.
We also remark that
the uniqueness statement in Theorem \ref{thm:main} could be sharpened, see
Section \ref{sect:proof_thm_main}.

In the exposition we focus on the DNLS equation but we remark that
the same approach is also applicable to slightly more general nonlinearities, cp. \cite{Tak99},
e.g.
$$\lambda_1 |u|^2 \partial_x u + \lambda_2 u^2 \partial_x \overline{u}+polynomial$$

We remark that the general strategy of proof of the tri-linear estimate is also applicable in the
non-periodic case \cite{Tak99}.

To illustrate the principle which allows us to gain the derivative
on the complex conjugate wave, let us consider three solutions $u_1,u_2,u_3$ of the
linear equation. Their Fourier transforms are supported on the parabola
$\{(\tau,\xi)\mid \tau+\xi^2=0\}$. The Fourier transform of the
interaction of the two linear waves $u_1,u_2$ at fixed frequencies $\xi_1,\xi_2$ with
$\partial_x\overline{u}_3$ is supported on $\{(\tau,\xi)\mid \tau+\xi^2=2(\xi-\xi_1)(\xi-\xi_2)\}$.
In the case where $\xi_1,\xi_2$ are small and the frequency $\xi_3$ is very
large, the frequency of the
resulting wave is also very large $\xi \sim \xi_3$. Hence its support is far away from the
parabola, or more precisely  $|\tau+\xi^2|\sim \xi_3^2$. This indicates that the Fourier
restriction norm method allows us gain a factor of order $\xi_3$ and everything reduces to terms
which can be treated by the $L^4$ Strichartz estimate \cite{Bo93}.
Moreover, we are able to control all other possible nonlinear interactions, as long
as $s\geq \frac{1}{2}$. 

The outline of the paper is as follows:
We conclude this Section with some general notation. In Section \ref{sect:gauge_trafo}
we introduce the gauge transform to link the DNLS with another derivative
nonlinear Schr\"odinger equation. After the introduction of useful function
spaces and linear estimates in Section \ref{sect:lin} we are concerned with
multi-linear estimates in Section
\ref{sect:multi}, which are applied in
Section \ref{sect:thm_gauge_eq_problem} to derive the sharp well-posedness
result for the gauge equivalent equation via the contraction mapping principle.
The proof of well-posedness for the DNLS is carried out in Section
\ref{sect:proof_thm_main}. Finally, the Appendix provides proofs of some
technical lemmata.

The author is indebted to Professor Herbert Koch, in particular
for helpful discussions about the gauge transform. Moreover, the author is
grateful to Axel Gr\"unrock and Martin Hadac for useful remarks.
\subsection*{Notation}
\noindent
Let $\mathcal{S}(\R)$ be the space of Schwartz functions.
We denote by $\test$ the space of functions $f: \R^2\to \C$
such that for all $(t,x) \in \R^2$
$$
f(t,x+2\pi)=f(t,2\pi),\; f(t, \cdot) \in C^{\infty}(\R), \; f(\cdot,x) \in
\mathcal{S}(\R)
$$

\noindent
We write $f_x=\partial_x f$ or $f_t=\partial_t f$ for partial
derivatives.

\noindent
Throughout this work $\chi \in C_0^{\infty}((-2,2))$ denotes a symmetric function
with $\chi\equiv 1$ in $[-1,1]$ and $\chi_T(t)=\chi(t/T)$.

\noindent
The Fourier transform with respect to the
spatial variable is defined by
$$
\mathcal{F}_x f (\xi)=(2\pi)^{-\frac{1}{2}}\int_0^{2\pi} e^{-i x\xi} f(x)\, dx \quad, \,\xi \in \Z
$$
and with respect to the time variable by
$$
\mathcal{F}_t f(\tau)=(2\pi)^{-\frac{1}{2}}\int_{\R} e^{-i t \tau}
f(t)\, dt \quad, \,\tau \in \R
$$
and $\mathcal{F}=\mathcal{F}_t\mathcal{F}_x$.

\noindent
For $1\leq p,q \leq \infty$ we use the notation
$$
\|f\|_{L^p_TL^q_x}:=\left\|t \mapsto \|x \mapsto f(t,x)\|_{L^q(\R)}
\right\|_{L^p([-T,T])}
$$
and if $T=\infty$ we write $\|f\|_{L^p_tL^q_x}$. Moreover,
$\|f\|_{L^p(\T)}=\|f\|_{L^p([0,2\pi])}$.

\noindent
We define the Sobolev spaces $H^s(\T)$ as the completion of the space of
$2\pi$-periodic $C^{\infty}$ functions $f$ with respect to the norm
\begin{equation*}
\|f\|^2_{H^s(\T)}:=\|f\|^2_{H^s}:= \sum_{\xi \in \Z} \lb \xi \rb^{2s}
  |\mathcal{F}_x f (\xi)|^2
\end{equation*}
where $\lb \xi
\rb=(1+|\xi|^2)^{\frac{1}{2}}$.

\noindent
The unitary group associated to $\partial_t - i\partial^2_{x}$ is defined via
$$\mathcal{F}_xW(t)u_0(\xi)=e^{-it\xi^2}\mathcal{F}_x
u_0(\xi)$$

\noindent
The Operator $J^s=J^s_x$ is defined by
$\mathcal{F}_xJ^sf(\xi)=\lb\xi\rb^s\mathcal{F}_xf(\xi)$.
We also use $J^s_t$ which is $J^s$ applied with respect to the $t$ variable.

\noindent
For $u \in C([-T,T],L^2(\T))$ we define $\mu(u)=\frac{1}{2\pi}\|u(0)\|^2_{L^2}$.

\noindent
For $\mu \in \R$ we define translations $\tau_{\mu}
u(t,x):=u(t,x+2\mu t)$ for $u \in C([-T,T],L^2(\T))$.

\section{The gauge transformation}\label{sect:gauge_trafo}
The Cauchy problem \eqref{eq:pdnls} is easily reduced to the case $\lambda=1$ by the transformation
$$
u(t,x) \mapsto \frac{1}{\sqrt{|\lambda|}} u(t,\sign(\lambda)x)
$$
From now on we only consider the case $\lambda=1$ without further comments.

Let $u \in C([-T,T],L^2(\T))$. We define the periodic primitive of $|u|^2-\frac{1}{2\pi}\|u(t)\|^2_{L^2}$ with
zero mean by
$$
\mathcal{I}(u)(t,x):=\frac{1}{2\pi}\int_0^{2\pi}\int_{\theta}^x |u(t,y)|^2
-\frac{1}{2\pi} \|u(t)\|_{L^2(\T)}^2 \, dy d\theta
$$
and
$v=e^{-i\mathcal{I}(u)}u(t,x)$. Now suppose that $u$ is a smooth solution to
\eqref{eq:pdnls} and let us derive an equation for $v$
\begin{align*}
v_t&=\exp(-i\mathcal{I}(u))(-i\mathcal{I}(u)_t u+u_t)\\
v_{xx}&=\exp(-i\mathcal{I}(u))(-\mathcal{I}(u)_x^2 u
-i\mathcal{I}(u)_x u_x-i(\mathcal{I}(u)_x u)_x+u_{xx})
\end{align*}
By the $L^2$
conservation law we have
$\|u(t)\|_{L^2(\T)}=\|u(0)\|_{L^2(\T)}$, see Appendix \ref{app_sect:cons_q}.
With 
$$
\mu:=\mu(u)=\frac{1}{2\pi}\|u(0)\|^2_{L^2(\T)}
$$
we have $$\mu(u)=\mu(v)=\frac{1}{2\pi}\|v(0)\|^2_{L^2(\T)}=\frac{1}{2\pi}\|v(t)\|^2_{L^2(\T)}$$
and $\mathcal{I}(u)_x(t,x)=|u(t,x)|^2-\mu$.
Therefore,
\begin{align}
&v_t-iv_{xx}\nonumber\\
&= e^{-i\mathcal{I}(u)}\left(u_t-iu_{xx}-(\mathcal{I}(u)_x u)_x
+i\mathcal{I}(u)_x^2 u-\mathcal{I}(u)_x u_x-i\mathcal{I}(u)_t u\right)\nonumber\\
&= e^{-i\mathcal{I}(u)}\left(\mu
  u_x+i(|u|^2-\mu)^2u-(|u|^2-\mu)u_x-i\mathcal{I}(u)_t u\right)
\label{eq:eq_v}
\end{align}
Moreover,
\begin{align*}
\diff{t}& \int_{\theta}^x |u(t,y)|^2-\mu \, dy=\int_{\theta}^x
u_t\overline{u}(t,y)+u\overline{u}_t(t,y)\, dy\\
=&\int_{\theta}^x
\left(iu_{xx}\overline{u}-i\overline{u}_{xx}u+\overline{u}(|u|^2u)_x+u(|u|^2\overline{u})_x\right)(t,y)\, dy
\end{align*}
Integration by parts yields
\begin{equation*}
 \int_{\theta}^xiu_{xx}\overline{u}(t,y)-i\overline{u}_{xx}u(t,y)\,
  dy=2\Imag(\overline{u}_{x}u)(t,x)-2\Imag(\overline{u}_{x}u)(t,\theta)
\end{equation*}
and
\begin{equation*}
\int_{\theta}^x
\overline{u}(|u|^2u)_x(t,y)+u(|u|^2\overline{u})_x(t,y)\, dy=\frac{3}{2}|u|^4(t,x)-\frac{3}{2}|u|^4(t,\theta)
\end{equation*}
which shows
\begin{align*}
\mathcal{I}(u)_t=&2\Imag(\overline{u}_{x}u)(t,x)+\frac{3}{2}|u|^4(t,x)\\
&-\frac{1}{2\pi}\int_0^{2\pi}
2\Imag(\overline{u}_{x}u)(t,\theta)+\frac{3}{2}|u|^4(t,\theta)\,d\theta
\end{align*}
Let us define
$$\phi(u)(t):=\frac{1}{2\pi}\int_0^{2\pi}
2\Imag(\overline{u}_{x}u)(t,\theta)+\frac{3}{2}|u|^4(t,\theta)\,d\theta$$
Plugging this into \eqref{eq:eq_v} we arrive at
\begin{align*}
v_t-iv_{xx}=e^{-i\mathcal{I}(u)}
\big(2\mu u_x-2i\mu|u|^2u+i\mu^2u-u^2\overline{u}_{x}-\tfrac{i}{2}|u|^4u+i\phi(u)u\big)
\end{align*}
Using $|u|=|v|$ as well as $u_x=e^{i\mathcal{I}(u)}v_x+i(|u|^2-\mu)u$ we get
\begin{align*}
v_t-iv_{xx}=&
2\mu v_x+2i\mu(|v|^2-\mu)v-2i\mu|v|^2v+i\mu^2v-v^2\overline{v}_{x}\\
&+i|v|^2(|v|^2-\mu)v-\tfrac{i}{2}|v|^4v+i\phi(u)v\\
=&2\mu v_x-v^2\overline{v}_{x}+\tfrac{i}{2}|v|^4v-i\mu|v|^2v+i(\phi(u)-\mu^2)v
\end{align*}
With
$$\psi(v)(t)=\frac{1}{2\pi}\int_0^{2\pi}
2\Imag(\overline{v}_{x}v)(t,\theta)-\frac{1}{2}|v|^4(t,\theta)\,d\theta
+\frac{1}{4\pi^2}\|v(0)\|^4_{L^2(\T)}
$$
we have $$\psi(v)=\phi(u)-\mu^2$$
and therefore obtain
\begin{equation*}
v_t-iv_{xx}-2\mu v_x
=-v^2\overline{v}_{x}+\tfrac{i}{2}|v|^4v-i\mu(v)|v|^2v+i\psi(v)v
\end{equation*}
Now, we use the transformation $w(t,x):=\tau_{-\mu}v(t,x):=v(t,x-2\mu
t)$ to cancel the linear term $2\mu v_x$ and arrive at
\begin{equation}\label{eq:gauge_equiv}
w_t-iw_{xx}
=-w^2\overline{w}_{x}+\tfrac{i}{2}|w|^4w-i\mu(w)|w|^2w+i\psi(w)w
\end{equation}
because $(\tau_{-\mu}v)_t(t,x)=v_t(t,x-2\mu
t)-2\mu v_x(t,x-2\mu t)$ and $\tau_{-\mu}$ commutes with partial differentiation
in $x$ as well as with $\psi$ and is an isometry in $L^2$.
The above calculation motivates the following definition.
\begin{definition}\label{def:gauge_trafo}
For $f \in L^2(\T)$ we define
$$\mathcal{G}(f)(x)=e^{-i\mathcal{I}(f)}f\,(x)
$$
where
$$
\mathcal{I}(f)(x):=\frac{1}{2\pi}\int_0^{2\pi}\int_{\theta}^x |f(y)|^2
-\frac{1}{2\pi} \|f\|_{L^2(\T)}^2 \, dy d\theta
$$
For $u \in C([-T,T],L^2(\T))$ we define
\begin{equation}\label{eq:gauge_trafo}
\mathcal{G}(u)(t,x):=\mathcal{G}(u(t))(x-2\mu(u) t)
\end{equation}
\end{definition}
\begin{remark}\label{rem:well_def}
The function $\mathcal{G}(f)$ is $2\pi$-periodic, since
$$|f(y)|^2
-\frac{1}{2\pi} \|f\|_{L^2(\T)}^2$$ has zero mean value and therefore
$$\int_0^{2\pi}\int_{\theta}^{x} |f(y)|^2
-\frac{1}{2\pi} \|f\|_{L^2(\T)}^2 \, dy d\theta$$
is $2\pi$-periodic.
\end{remark}

In the next Lemma, we summarize important properties of the nonlinear operator $\mathcal{G}$.
\begin{lemma}\label{lem:gauge_trafo_est}
For $s \geq 0$ the map
$$\mathcal{G}: C([-T,T],H^s(\T))\to
C([-T,T],H^s(\T))$$
is a homeomorphism.
Moreover, for $r>0$ there exists $c>0$, such that for
$$u,v \in B_{r,\mu}= \Big\{u \in C([-T,T],H^{s}(\T)) \mid \sup_{|t|\leq
    T}\|u(t)\|_{H^s(\T)}\leq r , \; \mu(u)=\mu \Big\}$$
the map $\mathcal{G}$ satisfies
\begin{equation}\label{eq:gauge_trafo_est}
\|\mathcal{G}(u)(t)-\mathcal{G}(v)(t)\|_{H^s(\T)}\leq c
\|u(t)-v(t)\|_{H^s(\T)} \quad , t \in [-T,T]
\end{equation}
for all $\mu \geq 0$.
The inverse map is given by
$$\mathcal{G}^{-1}(v)(t,x)=e^{i\mathcal{I}(\tau_{\mu(v)}v)}\tau_{\mu(v)} v(t,x)$$
and $\mathcal{G}^{-1}$ satisfies the same estimate \eqref{eq:gauge_trafo_est} on
subsets $B_{r,\mu}$.
Hence,
$\mathcal{G}$ is locally bi-Lipschitz on subsets with prescribed
$\|u(0)\|_{L^2}$.
\end{lemma}
\begin{proof}

We fix $s \geq 0$. There exists $c>0$, such that for $f,g,h \in H^s(\T)$
\begin{equation}\label{eq:exp_est}
\left\|(e^{\pm i \mathcal{I}(f)}-e^{\pm i
  \mathcal{I}(g)})h\right\|_{H^s}\leq c
e^{c\|f\|^2_{H^s}+c\|g\|^2_{H^s}}(\|f\|_{H^s}+\|g\|_{H^s})
\|f-g\|_{H^{s}}\|h\|_{H^s}
\end{equation}
This is proved in Appendix \ref{app_sect:exp_est}.
To show the Lipschitz estimate \eqref{eq:gauge_trafo_est} let $u,v \in B_{r,\mu}$. 
We observe that for fixed
$t$ a translation in $x$ is an isometric isomorphism on $H^s(\T)$. Using \eqref{eq:exp_est}
\begin{align*}
 &\|\mathcal{G}(u)(t)-\mathcal{G}(v)(t)\|_{H^s(\T)}\\
\leq&\left\|(e^{-i \mathcal{I}(u(t))}-e^{-i
  \mathcal{I}(v(t))})u(t)\right\|_{H^s}+\left\|(e^{-i
  \mathcal{I}(v(t))}-1)(u-v)(t)\right\|_{H^s}\\
&+\left\|(u-v)(t)\right\|_{H^s} \\
\leq &(2cre^{2cr^2}+cre^{cr^2}+1)\left\|(u-v)(t)\right\|_{H^s}
\end{align*}
which shows the Lipschitz continuity on $B_{r,\mu}$.

If $v=\mathcal{G}(u)$, then $\|v(0)\|_{L^2}=\|u(0)\|_{L^2}$ and therefore
$\mu(u)=\mu(v)$. Moreover, $|v(t,x+2\mu t)|=|u(t,x)|$ for a.e. $x$. Now, the
inversion formula is obvious and for
$\mathcal{G}^{-1}$ the Lipschitz estimate on subsets $B_{r,\mu}$
follows as above by replacing $-$ by $+$ in the exponential.

Now, the continuity of
$$
\mathcal{G},\mathcal{G}^{-1}: C([-T,T],H^s(\T)) \to C([-T,T],H^s(\T))
$$
on the whole space follows 
from the Lipschitz continuity of $\mathcal{G},\mathcal{G}^{-1}$ on subsets
$B_{r,\mu}$ and the continuity of the translations
$$
\tau_\mu: \R \to C([-T,T],H^s(\T)), \;\tau_\mu u (t,x)= u(t,x+2\mu t)
$$
together with the continuity of $\mu: C([-T,T],H^s(\T)) \to \R, \;\mu(u)=\frac{1}{2\pi}\|u(0)\|^2_{L^2}$.
\end{proof}

By repeating computation from the beginning of this section in reverse order, we prove
\begin{lemma}\label{lem:gauge_eq_prob}
Let $u,v \in C([-T,T],H^3(\T)) \cap C^1([-T,T],H^1(\T))$ and
$v=\mathcal{G}(u)$.
Then,
$u$ is a solution of
\begin{align*}
\partial_t u-i\partial^2_{x}u&=\partial_{x}(|u|^2u) \quad\text{in } (-T,T)\times \T\\
u(0)&=u_0
\end{align*}
if and only if $v$ is a solution of
\begin{align*}
\partial_t v-i\partial^2_{x}v
&=-v^2\partial_x\overline{v}+\tfrac{i}{2}|v|^4v-i\mu(v) |v|^2v+i\psi(v)v \quad\text{in } (-T,T)\times \T\\
v(0)&=\mathcal{G}(u_0)
\end{align*}
where
\begin{equation}\label{eq:psi}
\psi(v)(t)=\frac{1}{2\pi}\int_0^{2\pi}
2\Imag(\overline{v}_{x}v)(t,\theta)-\frac{1}{2}|v|^4(t,\theta)\,d\theta
+\mu(v)^2
\end{equation}
\end{lemma}

Finally, we show an estimate for $\psi$.
\begin{lemma}
Let $\psi$ be defined by \eqref{eq:psi}. Then,
\begin{equation}
  \label{eq:psi_est}
\begin{split}
|\psi(u)(t)-\psi(v)(t)|\leq&
c\big(1+\|u(t)\|_{H^{\frac{1}{2}}}
+\|v(t)\|_{H^{\frac{1}{2}}}\big)^3\|(u-v)(t)\|_{H^{\frac{1}{2}}}\\
+&2
(\|u(0)\|_{L^2}^3+\|v(0)\|_{L^2}^3)\|u(0)-v(0)\|_{L^2}
\end{split}
\end{equation}
\end{lemma}
\begin{proof}
We suppress the $t$ dependence and just write $u=u(t)$, $v=v(t)$.
\begin{equation*}
\left|\int_0^{2\pi}
(\Imag(\overline{u}_{x}u)-\Imag(\overline{v}_{x}v))(x)\,dx\right|\leq \left|\left(
u-v,u_{x}\right)_{L^2}\right|+ \left|\left(
v,u_{x}-v_{x}\right)_{L^2}\right|
\end{equation*}
Since $J_x^{\frac{1}{2}}$ is formally self-adjoint with respect to
$(\cdot,\cdot)_{L^2}$
we get
\begin{align*}
&\left|\left(
u-v,u_{x}\right)_{L^2}
\right|
+\left|\left(
v,u_{x}-v_{x}\right)_{L^2}\right|\\
=&\left|\left(
J_x^{\frac{1}{2}}(u-v),J_x^{-\frac{1}{2}}\partial_x u\right)_{L^2}\right|+\left|\left(
J_x^{\frac{1}{2}}v,J_x^{-\frac{1}{2}}\partial_x(u-v)\right)_{L^2}\right|\\
\leq &c (\|u\|_{H^\frac{1}{2}}+\|v\|_{H^\frac{1}{2}})\|u-v\|_{H^{\frac{1}{2}}} 
\end{align*}
Moreover,
\begin{align*}
\left|\int_0^{2\pi}
(|u|^4-|v|^4)(x)\,dx\right|\leq &\int_0^{2\pi}
\left||u|-|v|\right|(|u|^3+|u|^2|v|+|u||v|^2+|v|^3)(x)\,dx
\\
\leq &2 (\|u\|_{L^6}^3+\|v\|_{L^6}^3)\|u-v\|_{L^2}
\end{align*}
Finally,
$$\left|\|u(0)\|^4_{L^2}-\|v(0)\|^4_{L^2}\right|\leq 2
(\|u(0)\|_{L^2}^3+\|v(0)\|_{L^2}^3)\|u(0)-v(0)\|_{L^2}$$
These three estimates together with the Sobolev embedding $H^{\frac{1}{3}}\hookrightarrow L^6$ prove \eqref{eq:psi_est}.
\end{proof}
\section{Definition of the spaces and linear estimates}\label{sect:lin}
The following spaces are well-known from \cite{Bo93,G96,GTV97,CKSTT03,G00}.
\begin{definition}\label{def:spaces}
Let $s,b \in \R$. The Bourgain space
$X_{s,b}$ associated to the Schr\"odinger
operator $\partial_t - i\partial^2_{x}$ is defined
as the completion of the space $\test$ with respect to the norm
\begin{equation}
  \label{eq:X_s,b}
  \|f\|^2_{X_{s,b}}:= \sum_{\xi \in \Z}\int_\R  \lb \xi \rb^{2s} \lb \tau +
  \xi^2\rb^{2b}  |\mathcal{F} f (\tau,\xi)|^2\, d\tau
\end{equation}
Similarly we define $X^-_{s,b}$ by replacing $\lb \tau +
  \xi^2\rb$ with $\lb \tau -\xi^2\rb$.

Moreover,
we define $Y_{s,b}$ as the completion of the space $\test$ with respect to
\begin{equation}
  \label{eq:Y_s}
  \|f\|^2_{Y_{s,b}}:= \sum_{\xi \in \Z} \left(\int_\R \lb \tau +
  \xi^2\rb^{b} \lb \xi \rb^{s} |\mathcal{F} f (\tau,\xi)|\, d\tau\right)^2
\end{equation}
and the space $Z_{s}:=X_{s,\frac{1}{2}} \cap Y_{s,0}$ with norm
\begin{equation}
  \label{eq:Z}
  \|u\|_{Z_{s}}:= \|u\|_{X_{s,\frac{1}{2}}} + \|u\|_{Y_{s,0}}
\end{equation}

For $T>0$ we define the restriction norm space 
$$Z^{T}_{s}:=\{u|_{[-T,T]} \mid u \in Z_{s} \}$$
with norm
$$\|u\|_{Z^{T}_{s}}=\inf\{\|\widetilde{u}\|_{Z_{s}} \mid u=
\widetilde{u}|_{[-T,T]}, \, \widetilde{u}\in Z_{s}\},$$

Finally, we define the metric space
$$
\mathcal{X}_s^T:=\left\{u \in C([-T,T],H^s(\T)) \mid \mathcal{G}(u) \in Z^{T}_s \right\}
$$
with metric
$d_{\mathcal{X}_s^T}(u,v)=\|\mathcal{G}(u)-\mathcal{G}(v)\|_{Z^{T}_s}$.
\end{definition}
\begin{remark}\label{rem:spaces}
\begin{enumerate}
\item For complex conjugation we observe $\|\overline{u}\|_{X_{s,b}}=\|u\|_{X^-_{s,b}}$.
\item The metric space $\mathcal{X}_s^T$ is complete.
\end{enumerate}
\end{remark}
Now, we start with frequently used embedding theorems.
\begin{lemma}\label{lem:sob}
\begin{align}
\text{If }\;2 \leq p<\infty, b\geq
\frac{1}{2}-\frac{1}{p}\,:
\|u\|_{L^p_tH^s}\leq & c\|u\|_{X_{s,b}}  \label{eq:sob}\\
\text{If }\;2 \leq p,q<\infty, b\geq
\frac{1}{2}-\frac{1}{p}\;,\, s\geq
\frac{1}{2}-\frac{1}{q}\,:
\|u\|_{L^p_tL^q_x}\leq& c\|u\|_{X_{s,b}}  \label{eq:sob2}\\
\text{If }\;1 < p \leq 2\;,\, b\leq
\frac{1}{2}-\frac{1}{p}\,: \|u\|_{X_{s,b}} \leq &c
\|u\|_{L^p_tH^s}\label{eq:dual_sob}
\end{align}
We may replace $X_{s,b}$ by $X^-_{s,b}$.
Moreover,
\begin{align}
\|u\|_{C(\R,H^s(\T))} \leq & c\|u\|_{Z_s}\quad  , s\in \R\label{eq:linfty_emb}\\
\|u\|_{Y_{s,b_1}}\leq & c\|u\|_{X_{s,b_2}}\quad  , b_2>b_1+\frac{1}{2}\label{eq:x_in_y}
\end{align}
\end{lemma}

\begin{proof}
We consider $v=W(-t)J_x^su$ for $u\in \test$. Then, by
Minkowski's and Sobolev's inequality
$$\|u\|_{L^p_tH^s_x}=\|v\|_{L^p_tL^2_x}\leq \|v\|_{L^2_xL^p_t}\leq c\|J^b_tv\|_{L^2_xL^2_t}=
c\|u\|_{X_{s,b}}$$
and the claim \eqref{eq:sob} follows. Combining this with another application of
Sobolev's inequality in the space variable $\|v(t)\|_{L^q_x}\leq c\|J^s_x v(t)\|_{L^2_x}$ gives \eqref{eq:sob2}.
Estimate \eqref{eq:dual_sob} follows by
duality from \eqref{eq:sob}. The estimates for $X^-_{s,b}$ follow from the
invariance of $L^p_tH_x^s$ and $L^p_tL^q_x$ under complex conjugation.
To prove \eqref{eq:linfty_emb} it suffices to prove an estimate for the $\sup$
norm for $u \in \test$ by density. We write for $t\in \R$
$$
\mathcal{F}_x u(t,\xi)=c\int_\R
e^{it\tau}\mathcal{F} u(\tau,\xi)\, d\tau
$$
by the Fourier inversion formula. This yields
$$\|u(t)\|_{H^s}=c\left\|\int_\R
e^{it\tau} \lb\xi\rb^s \mathcal{F} u(\tau,\xi)\, d\tau\right\|_{L^2_\xi}\leq c\|
\lb\xi\rb^s \mathcal{F} u(\tau,\xi)\|_{L^2_{\xi} L^1_{\tau}}$$
Now we take the supremum with respect to $t$.
The last estimate follows from the Cauchy-Schwarz inequality in $\tau$:
\begin{align*}
\|u\|^2_{Y_{s,b_1}}=&\sum_{\xi \in \Z} \left(\int_\R \lb \tau +
  \xi^2\rb^{b_1} \lb \xi \rb^{s} |\mathcal{F} f (\tau,\xi)|\, d\tau\right)^2\\
 \leq &\sum_{\xi \in \Z} \int_\R \lb \tau +
  \xi^2\rb^{2b_1-2b_2}\, d\tau \, \int_\R \lb \tau +
  \xi^2\rb^{2b_2} \lb \xi \rb^{2s} |\mathcal{F} f (\tau,\xi)|^2\, d\tau
\end{align*}
Since by assumption $2b_1-2b_2<-1$, there exists $c>0$, such that for all $\xi$
$$\int_\R \lb \tau +
  \xi^2\rb^{2b_1-2b_2}\, d\tau \leq c
$$ 
which finishes the proof.
\end{proof}

\begin{lemma}\label{lem:strichartz}
For $-b',b>\frac{3}{8}$ there exists $c>0$, such that
  \begin{equation}
    \label{eq:strichartz}
    \|u\|_{L_{tx}^4} \leq c \|u\|_{X_{0,b}}
  \end{equation} 

\begin{equation}
    \label{eq:strichartz_dual}
     \|u\|_{X_{0,b'}} \leq c  \|u\|_{L_{tx}^{4/3}}
  \end{equation} 
In these estimates we may also replace $X_{0,b}$ by $X^-_{0,b}$.
\end{lemma}
\begin{proof}
The estimate \eqref{eq:strichartz} is essentially Bourgain's $L^4$ Strichartz estimate
\cite{Bo93}, but in a version which is global in time.
This can be found in \cite{G00}, Lemma 2.1.
Using duality this also shows \eqref{eq:strichartz_dual}.
That the estimates hold both for $X_{s,b}$ and $X^-_{s,b}$ results from the
invariance of $L^p_tL^q_x$ spaces under complex conjugation.
\end{proof}
We summarize the behavior of the $X_{s,b},Y_{s,0}$ norms under
multiplication with cutoffs in time. This will be frequently used in the
sequel without further remarks.
\begin{lemma}\label{lem:cutoffs}
Let $s\in \R$ and $0<T \leq 1$.
There exists $c>0$, such that
\begin{equation*}
\|\chi_T u\|_{Y_{s,0}} \leq c \|u\|_{Y_{s,0}}
\end{equation*}
Moreover, for $0 \leq b_1<b_2<\frac{1}{2}$ or $-\frac{1}{2}<b_1<b_2 \leq 0$ there exists $c>0$, such that
\begin{equation*}
\|\chi_T u\|_{X_{s,b_1}} \leq cT^{b_2-b_1}\|u\|_{X_{s,b_2}}
\end{equation*}
and for any $\delta>0$ there exists $c>0$, such that
\begin{equation*}
\|\chi_T u\|_{X_{s,\frac{1}{2}}} \leq c T^{-\delta}\| u\|_{X_{s,\frac{1}{2}}}
\end{equation*}
\end{lemma}
\begin{proof}
The first estimate follows from Young's inequality in $\tau$: For fixed $\xi$ we have
\begin{align*}
\|\mathcal{F}(\chi_T u)(\tau,\xi)\|_{L^1_\tau}=&c\|\int_\R \mathcal{F}_t\chi_T(\tau-\tau_1)
\mathcal{F}u(\tau_1,\xi)\, d\tau_1\|_{L^1_\tau}\\
\leq & c\|\mathcal{F}_t\chi_T\|_{L^1_\tau}
\|\mathcal{F}u(\tau,\xi)\|_{L^1_\tau}
\end{align*}
Because $\|\mathcal{F}_t\chi_T\|_{L^1_\tau}=\|\mathcal{F}_t\chi\|_{L^1_\tau}$
the estimate follows by taking the $L^2_\xi$ norms on both sides.
The second estimate is proved in \cite{G00}, Lemma 1.2 and for the third
estimate we refer to the proof of \cite{GTV97}, Lemma 2.5 and the subsequent
remark, which remain true in the periodic setting.
\end{proof}
The next Lemma contains the classical (cp. \cite{Bo93,G96,GTV97,CKSTT03})
estimates for the linear homogeneous and inhomogeneous problem.
\begin{lemma}\label{lem:lin_est}
Let $s \in \R$. There exists $c>0$, such that for all $u_0 \in H^s(\T)$
\begin{equation}\label{eq:lin_hom}
\left\|\chi W(t)u_0 \right\|_{Z_s} \leq c \|u_0\|_{H^s}
\end{equation}
and for all $f \in \test$ with $\supp(f) \subset \{(t,x)\mid |t|\leq 2\}$
\begin{equation}\label{eq:lin_inhom}
\left\|\chi \int_0^t W(t-t') f(t')\, dt'\right\|_{Z_s} \leq c \|f\|_{Y_{s,-1}}+ c \|f\|_{X_{s,-\frac{1}{2}}}
\end{equation}
\end{lemma}
\begin{proof}
It suffices to consider smooth $u_0$. Let us write
$$
\mathcal{F}(\chi W(\cdot)u_0)(\tau,\xi)=\mathcal{F}_t \chi
(\tau+\xi^2)\mathcal{F}_x u_0(\xi)
$$
Then, because $\mathcal{F}_t \chi$ is a Schwartz function
the estimate \eqref{eq:lin_hom} follows.
Now we turn to the estimate \eqref{eq:lin_inhom} for the linear inhomogeneous
equation and we follow the argumentation from \cite{CKSTT03}, Lemma 3.1.
We have
$$
\chi(t) \int_0^t W(t-t') f(t')\, dt'=I(t)+J(t)
$$
with
\begin{align*}
I(t)
=&\frac{1}{2}\chi(t)W(t)\int_\R \varphi(t') W(-t') f(t')\, dt'\\
J(t)=&\frac{1}{2}\chi(t)\int_\R \varphi(t-t')W(t-t') f(t')\, dt'
\end{align*}
and $\varphi(t')=\chi(t'/10)\sign(t')$.
Moreover,
\begin{equation}
  \label{eq:varphi}
  |\mathcal{F}_t \varphi(\tau)| \leq c \lb\tau\rb^{-1}
\end{equation}
Now, by estimate \eqref{eq:lin_hom}
$$
\|I\|_{Z_s} \leq c \left\|\int_\R \varphi(t') W(-t') f(t')\, dt'\right\|_{H^s(\T)}
$$
and by Parseval's equality
$$
\mathcal{F}_x \left(\int_\R \varphi(t') W(-t') f(t')\, dt'\right)(\xi) =\int_\R
\overline{\mathcal{F}_t\varphi}(\tau+\xi^2) \mathcal{F}f(\tau,\xi) \, d\tau
$$
which implies
$$\left\|\int_\R \varphi(t') W(-t') f(t')\, dt'\right\|_{H^s(\T)} \leq c \|f\|_{Y_{s,-1}}
$$
by \eqref{eq:varphi}. To show the estimate for $J$ we first apply Lemma \ref{lem:cutoffs} with $T=1$

$$\|J\|_{Z_s}\leq c\left\|\int_\R \varphi(t-t')W(t-t') f(t')\, dt'\right\|_{Z_s}
$$
and observe
$$
\mathcal{F}\left(\int_\R \varphi(t-t')W(t-t') f(t')\, dt'\right) (\tau,\xi)
= \mathcal{F}_t\varphi(\tau+\xi^2) \mathcal{F}f(\tau,\xi)
$$
The claim then follows from \eqref{eq:varphi}.
\end{proof}
\section{Multi-linear estimates}\label{sect:multi}
We start with an elementary bound for the multi-linear multiplier in the spirit of \cite{CKSTT03}.
\begin{lemma}\label{lem:multiplier_bound1}
Let $\tau_j \in \R$ and $\xi_j \in \Z$, $j=1,2,3$
and define $$M(\tau_1,\tau_2,\tau_3,\xi_1,\xi_2,\xi_3)
=\frac{\lb\xi\rb^{\frac{1}{2}} i \xi_3}{
\lb\tau+\xi^2\rb^{\frac{1}{2}}
\lb\tau_1+\xi_1^2\rb^{\frac{1}{2}}\lb\tau_2+\xi_2^2\rb^{\frac{1}{2}}\lb\tau_3-\xi_3^2\rb^{\frac{1}{2}}
\prod_{j=1}^3\lb\xi_j\rb^{\frac{1}{2}}}$$
where $(\tau,\xi)=(\tau_1,\xi_1)+(\tau_2,\xi_2)+(\tau_3,\xi_3)$. Moreover let
\begin{align*}
M_0(\tau_1,\tau_2,\tau_3,\xi_1,\xi_2,\xi_3)=&\frac{\chi_{A_0}}
{\lb\xi_1\rb^{\frac{1}{2}}\lb\xi_2\rb^{\frac{1}{2}}
\lb\tau_1+\xi_1^2\rb^{\frac{1}{2}}\lb\tau_2+\xi_2^2\rb^{\frac{1}{2}}\lb\tau_3-\xi_3^2\rb^{\frac{1}{2}}}\\
M_j(\tau_1,\tau_2,\tau_3,\xi_1,\xi_2,\xi_3)=&\frac{\chi_{A_j}}
{\lb\xi_1\rb^{\frac{1}{2}}\lb\xi_2\rb^{\frac{1}{2}}\lb\tau+\xi^2\rb^{\frac{1}{2}}
\lb\tau_k+\xi_k^2\rb^{\frac{1}{2}}\lb\tau_3-\xi_3^2\rb^{\frac{1}{2}}}
\end{align*}
for $j,k \in \{1,2\}$ and $k\not= j$, as well as
\begin{align*}
M_3(\tau_1,\tau_2,\tau_3,\xi_1,\xi_2,\xi_3)=&\frac{\chi_{A_3}}
{\lb\xi_1\rb^{\frac{1}{2}}\lb\xi_2\rb^{\frac{1}{2}}\lb\tau+\xi^2\rb^{\frac{1}{2}}
\lb\tau_1+\xi_1^2\rb^{\frac{1}{2}}\lb\tau_2+\xi_2^2\rb^{\frac{1}{2}}}\\
N(\tau_1,\tau_2,\tau_3,\xi_1,\xi_2,\xi_3)
=&\frac{1}{\lb\tau+\xi^2\rb^{\frac{1}{2}}\lb\tau_1+\xi_1^2\rb^{\frac{1}{2}}
\lb\tau_2+\xi_2^2\rb^{\frac{1}{2}}\lb\tau_3-\xi_3^2\rb^{\frac{1}{2}}}
\end{align*}
where $\chi_{A_0}$ is the characteristic function of the subregion $A_0 \subset \R^3 \times \Z^3$ where
$$\lb\tau+\xi^2\rb \geq
\lb\tau_1+\xi_1^2\rb,\lb\tau_2+\xi_2^2\rb,\lb\tau_3-\xi_3^2\rb$$
and the subregions $A_j$ for $j=1,2,3$  are defined analogously.
Then, the estimate
\begin{equation}\label{eq:multiplier_bound1}
|M|\leq 16( \sum_{j=0}^3 M_j+N)
\end{equation}
holds true.
\end{lemma}
\begin{proof}
The key for the proof will be the observation that
\begin{equation}
  \label{eq:resonance}
\tau+\xi^2- \left(\tau_1+\xi_1^2+\tau_2+\xi_2^2+\tau_3-\xi_3^2\right)=2(\xi-\xi_1)(\xi-\xi_2)
\end{equation}
which implies
\begin{align}
&\lb(\xi-\xi_1)(\xi-\xi_2)\rb^{\frac{1}{2}} \leq 
 \lb\tau+\xi^2\rb^{\frac{1}{2}}+\lb\tau_1+\xi_1^2\rb^{\frac{1}{2}}
+\lb\tau_2+\xi_2^2\rb^{\frac{1}{2}}+\lb\tau_3-\xi_3^2\rb^{\frac{1}{2}}\nonumber\\
 &\leq 4\left(\chi_{A_0}\lb\tau+\xi^2\rb^{\frac{1}{2}}+ \chi_{A_1}\lb\tau_1+\xi_1^2\rb^{\frac{1}{2}}
+ \chi_{A_2}\lb\tau_2+\xi_2^2\rb^{\frac{1}{2}}+ \chi_{A_3}\lb\tau_3-\xi_3^2\rb^{\frac{1}{2}}\right)
\label{eq:res_est}
\end{align}
We distinguish 4 different cases.
\begin{enumerate}
\item\label{case:1} $|\xi| > 2 |\xi_1|$ and $|\xi| > 2 |\xi_2|$:
In this case $|\xi_3|\leq 2 |\xi|$ and $4 \lb(\xi-\xi_1)(\xi-\xi_2)\rb \geq \lb
\xi \rb^2$ and we conclude $|M|\leq 16 \sum_{j=0}^3 M_j$.
\item\label{case:2} $|\xi|\leq 2 |\xi_1|$ and $|\xi|\leq 2 |\xi_2|$:
In this case we have $|\xi_3|\leq 4 \max\{|\xi_1|,|\xi_2|\}$ and $|\xi|\leq 2
\min\{|\xi_1|,|\xi_2|\}$, which shows $|M|\leq 4 N$.
\item\label{case:3} $|\xi|> 2 |\xi_1|$ and $|\xi|\leq 2 |\xi_2|$:
We have $$|\xi| \leq 2|\xi-\xi_1| \text{ and }
|\xi||\xi-\xi_2|\leq 2 \lb(\xi-\xi_1)(\xi-\xi_2)\rb$$
Due to the fact that
$$|\xi_3|^{\frac{1}{2}}\leq
|\xi-\xi_2|^{\frac{1}{2}}+|\xi_1|^{\frac{1}{2}}$$
we have
$$\lb\xi\rb^{\frac{1}{2}}|\xi_3|^{\frac{1}{2}}\leq 2|\xi|^{\frac{1}{2}}|\xi_3|^{\frac{1}{2}}
\leq 2|\xi|^{\frac{1}{2}}|\xi-\xi_2|^{\frac{1}{2}}+2|\xi|^{\frac{1}{2}}|\xi_1|^{\frac{1}{2}}
$$
The first term is bounded by $4\lb(\xi-\xi_1)(\xi-\xi_2)\rb^{\frac{1}{2}}$ and we apply
\eqref{eq:res_est}. The second term is bounded
by
$4\lb\xi_2\rb^{\frac{1}{2}}\lb\xi_1\rb^{\frac{1}{2}}$,
which proves $|M|\leq 16 \sum_{j=0}^3 M_j+4 N$.
\item\label{case:4} $|\xi|\leq 2 |\xi_1|$ and $|\xi|> 2 |\xi_2|$:
By the symmetry of $M$ in $\xi_1,\xi_2$ we find the same estimate as in case \ref{case:3}.
\end{enumerate}
\end{proof}

In the sequel we will use the abbreviations
\begin{align*}
&\int_\ast\sum_{\ast} \prod_{j=1}^k
f_j(\tau_j,\xi_j):=\int\limits_{\tau=\tau_1+\ldots+\tau_k}\sum\limits_{\xi=\xi_1+\ldots+\xi_k}
\prod_{j=1}^k f_j(\tau_j,\xi_j)\\
:=
&\int\limits_{\R^{k-1}}\sum\limits_{(\xi_1,\ldots,\xi_{k-1}) \in \Z^{k-1}}
\prod_{j=1}^{k-1}
f_j(\tau_j,\xi_j)f_k(\tau-\sum_{j=1}^{k-1}\tau_j,\xi-\sum_{j=1}^{k-1}\xi_j)
d\tau_1 \ldots d\tau_{k-1}
\end{align*} which is nothing else but the convolution $f_1 \ast \ldots \ast f_k (\tau,\xi)$.
\begin{theorem}\label{thm:trilinear_est1}
There exists $c,\eps>0$, such that for $T \in (0,1]$ and $u_j \in
\test $ with $\supp(u_j) \subset \{(t,x)\mid
|t|\leq T\}$, $j=1,2,3$, we have 
\begin{equation}
  \label{eq:trilinear_est1}
  \left\|u_1 u_2 \partial_x u_3\right\|_{X_{\frac{1}{2},-\frac{1}{2}}}
\leq c
T^{\eps}\|u_1\|_{X_{\frac{1}{2},\frac{1}{2}}}\|u_2\|_{X_{\frac{1}{2},\frac{1}{2}}}
\|u_3\|_{X^-_{\frac{1}{2},\frac{1}{2}}}
\end{equation}
\end{theorem}

\begin{proof}
We define
$$
f_j(\tau_j,\xi_j)=\lb\tau_j+\xi_j^2\rb^{\frac{1}{2}}\lb\xi_j\rb^{\frac{1}{2}}\mathcal{F}u_j(\tau_j,\xi_j)
$$
for $j=1,2$ and
$$
f_3(\tau_3,\xi_3)=\lb\tau_3-\xi_3^2\rb^{\frac{1}{2}}\lb\xi_3\rb^{\frac{1}{2}}\mathcal{F}u_3(\tau_3,\xi_3)
$$
With the Fourier multiplier $M$ defined in Lemma
\ref{lem:multiplier_bound1}
we rewrite the left hand side as
\begin{equation*}
\left\|u_1 u_2 \partial_x u_3\right\|_{X_{\frac{1}{2},-\frac{1}{2}}}=\left\|
\int_\ast
\sum_\ast M(\tau_1,\tau_2,\tau_3,\xi_1,\xi_2,\xi_3)
\prod_{j=1}^3 f_j(\tau_j,\xi_j)
\right\|_{l^2_{\xi}L^2_{\tau}}
\end{equation*}
By an application of the triangle inequality we may assume $f_j,
\mathcal{F}u_j \geq 0$ and $\|u_j\|_{X_{s,b}}=\|\chi_T u_j\|_{X_{s,b}}$.
By the point-wise bound
\eqref{eq:multiplier_bound1} on $|M|$ the left hand side is
bounded by the sum over the corresponding terms with $M$ replaced by
$M_0,M_1,M_2,M_3$ and $N$, respectively.\\
\noindent\emph{Estimate for $M_0$:}
\begin{align*}
&\left\|\int_\ast
\sum_\ast M_0(\tau_1,\tau_2,\tau_3,\xi_1,\xi_2,\xi_3)
\prod_{j=1}^3 f_j(\tau_j,\xi_j)
\right\|_{l^2_{\xi}L^2_{\tau}}
=\|u_1u_2J^{\frac{1}{2}}u_3\|_{L^2_tL^2_x}=:m_0
\end{align*}
Using H\"older, we get
\begin{align*}
m_0 \leq\|u_1\|_{L^8_tL^8_x}\|u_2\|_{L^8_tL^8_x}\|J^{\frac{1}{2}}u_3\|_{L^4_tL^4_x} \leq
 c\|u_1\|_{X_{\frac{3}{8},\frac{3}{8}}}\|u_2\|_{X_{\frac{3}{8},\frac{3}{8}}}
\|u_3\|_{X^-_{\frac{1}{2},\frac{1}{2}}}
\end{align*}
where we used Sobolev's inequality w.r.t. space and time variables on $u_1,u_2$
as well as the $L^4$ Strichartz inequality on $J^{\frac{1}{2}}u_3$. By the
localization in time, see Lemma \ref{lem:cutoffs}
$$
m_0 \leq
cT^{\eps}\|u_1\|_{X_{\frac{1}{2},\frac{1}{2}}}\|u_2\|_{X_{\frac{1}{2},
\frac{1}{2}}}\|u_3\|_{X^-_{\frac{1}{2},\frac{1}{2}}}
$$
\noindent\emph{Estimate for $M_1$:}
\begin{align}
&\left\|\int_\ast
\sum_\ast M_1(\tau_1,\tau_2,\tau_3,\xi_1,\xi_2,\xi_3)
\prod_{j=1}^3 f_j(\tau_j,\xi_j)
\right\|_{l^2_{\xi}L^2_{\tau}}\nonumber\\
=&
\|J^{-\frac{1}{2}}\mathcal{F}^{-1}f_1u_2J^{\frac{1}{2}}u_3\|_{X_{0,-\frac{1}{2}}}
\leq \|J^{-\frac{1}{2}}\mathcal{F}^{-1}f_1u_2J^{\frac{1}{2}}u_3\|_{X_{0,-\frac{3}{8}}}
=:m_1\label{eq:m_1}
\end{align}
Then, by Sobolev in time
\begin{align*}
m_1 \leq& c
\|J^{-\frac{1}{2}}\mathcal{F}^{-1}f_1u_2J^{\frac{1}{2}}u_3\|_{L^{8/7}_tL^{2}_x}\\
\leq &
 c\|J^{-\frac{1}{2}}\mathcal{F}^{-1}f_1\|_{L^2_tL^{8}_x}
\|u_2J^{\frac{1}{2}}u_3\|_{L^{8/3}_tL^{8/3}_x}\\
\leq &c\|J^{-\frac{1}{2}}\mathcal{F}^{-1}f_1\|_{L^2_tL^{8}_x}\|u_2\|_{L^{8}_tL^{8}_x}
\|J^{\frac{1}{2}}u_3\|_{L^{4}_tL^{4}_x}
\end{align*}
Now we use the Sobolev inequality on the first two factors
as well as the $L^4$ Strichartz inequality on $J^{\frac{1}{2}}u_3$ and obtain
\begin{equation}
m_1 \leq
cT^{\eps}\|u_1\|_{X_{\frac{1}{2},\frac{1}{2}}}\|u_2\|_{X_{\frac{1}{2},
\frac{1}{2}}}\|u_3\|_{X^-_{\frac{1}{2},\frac{1}{2}}}\label{eq:bound_m_1}
\end{equation}
\noindent\emph{Estimate for $M_2$:}
\begin{align}
&\left\|\int_\ast
\sum_\ast M_2(\tau_1,\tau_2,\tau_3,\xi_1,\xi_2,\xi_3)
\prod_{j=1}^3 f_j(\tau_j,\xi_j)
\right\|_{l^2_{\xi}L^2_{\tau}}\nonumber\\
=&
\|u_1J^{-\frac{1}{2}}\mathcal{F}^{-1}f_2J^{\frac{1}{2}}u_3\|_{X_{0,-\frac{1}{2}}}\leq
\|u_1J^{-\frac{1}{2}}\mathcal{F}^{-1}f_2J^{\frac{1}{2}}u_3\|_{X_{0,-\frac{3}{8}}}=:m_2\label{eq:m_2}
\end{align}
As for $m_1$, by exchanging the roles of the first two factors we obtain
\begin{equation}
m_2 \leq
cT^{\eps}\|u_1\|_{X_{\frac{1}{2},\frac{1}{2}}}\|u_2\|_{X_{\frac{1}{2},
\frac{1}{2}}}\|u_3\|_{X^-_{\frac{1}{2},\frac{1}{2}}}\label{eq:bound_m_2}
\end{equation}
\noindent\emph{Estimate for $M_3$:}
\begin{align}
&\left\|\int_\ast
\sum_\ast M_3(\tau_1,\tau_2,\tau_3,\xi_1,\xi_2,\xi_3)
\prod_{j=1}^3 f_j(\tau_j,\xi_j)
\right\|_{l^2_{\xi}L^2_{\tau}}\nonumber\\
=&\|u_1u_2\mathcal{F}^{-1}f_3\|_{X_{0,-\frac{1}{2}}}
\leq \|u_1u_2\mathcal{F}^{-1}f_3\|_{X_{0,-\frac{7}{16}}}=:m_3
\label{eq:m_3}
\end{align}
We apply dual Strichartz' \eqref{eq:strichartz_dual}, H\"older's and Sobolev's
inequality to conclude
\begin{align}
  m_3 \leq& c \|u_1u_2\mathcal{F}^{-1}f_3\|_{L^{4/3}_tL^{4/3}_x}\nonumber\\
\leq& c\|u_1\|_{L^8_tL^8_x}\|u_2\|_{L^8_tL^8_x}\|f_3\|_{L^2_tL^2_x}\nonumber\\
\leq &cT^{\eps} \|u_1\|_{X_{\frac{1}{2},\frac{1}{2}}}\|u_2\|_{X_{\frac{1}{2},
\frac{1}{2}}}\|u_3\|_{X^-_{\frac{1}{2},\frac{1}{2}}}\label{eq:bound_m_3}
\end{align}
\noindent\emph{Estimate for $N$:}
\begin{align}
&\left\|\int_\ast
\sum_\ast N(\tau_1,\tau_2,\tau_3,\xi_1,\xi_2,\xi_3)
\prod_{j=1}^3 f_j(\tau_j,\xi_j)
\right\|_{l^2_{\xi}L^2_{\tau}}\nonumber\\
=&\|J^{\frac{1}{2}}u_1J^{\frac{1}{2}}u_2J^{\frac{1}{2}}u_3\|_{X_{0,-\frac{1}{2}}}
\leq \|J^{\frac{1}{2}}u_1J^{\frac{1}{2}}u_2J^{\frac{1}{2}}u_3\|_{X_{0,-\frac{7}{16}}}=:n\label{eq:n}
\end{align}
Strichartz inequalities \eqref{eq:strichartz} and \eqref{eq:strichartz_dual}
yield
\begin{align}
  n\leq & c
\|J^{\frac{1}{2}}u_1J^{\frac{1}{2}}u_2J^{\frac{1}{2}}u_3\|_{L^{4/3}_tL^{4/3}_x}\nonumber\\
\leq & c
 \|J^{\frac{1}{2}}u_1\|_{L^{4}_tL^{4}_x}\|J^{\frac{1}{2}}u_2\|_{L^{4}_tL^{4}_x}
\|J^{\frac{1}{2}}u_3\|_{L^{4}_tL^{4}_x}\nonumber\\
\leq & cT^{\eps}\|u_1\|_{X_{\frac{1}{2},\frac{1}{2}}}\|u_2\|_{X_{\frac{1}{2},
\frac{1}{2}}}\|u_3\|_{X^-_{\frac{1}{2},\frac{1}{2}}}\label{eq:bound_n}
\end{align}
and the proof is complete.
\end{proof}
\begin{lemma}\label{lem:multiplier_bound2} We use the notation from Lemma
  \ref{lem:multiplier_bound1} and
define $$\tilde{M}(\tau_1,\tau_2,\tau_3,\xi_1,\xi_2,\xi_3)
:=\frac{M(\tau_1,\tau_2,\tau_3,\xi_1,\xi_2,\xi_3)}{\lb\tau+\xi^2\rb^{\frac{1}{2}}}$$
Then,
\begin{equation}\label{eq:multiplier_bound2}
|\tilde{M}|\leq 128(\sum_{j=0}^3 \tilde{M}_j+\tilde{N})
\end{equation}
where for $\delta \in (0,\tfrac{1}{6})$
$$
\tilde{M}_0
=\frac{\chi_{A_0}}{\lb\xi\rb^{\frac{1}{2}-3\delta}\lb\xi_1\rb^{\frac{1}{2}}
\lb\xi_2\rb^{\frac{1}{2}}\lb\xi_3\rb^{\frac{1}{2}-3\delta}
\lb\tau_1+\xi_1^2\rb^{\frac{1}{2}+\delta}\lb\tau_2+\xi_2^2\rb^{\frac{1}{2}+\delta}
\lb\tau_3-\xi_3^2\rb^{\frac{1}{2}+\delta}}
$$
and
$$
\tilde{M}_j=\lb\tau+\xi^2\rb^{-\frac{1}{2}}M_j
\;,j \in \{1,2,3\}\quad,\; \tilde{N}:=\lb\tau+\xi^2\rb^{-\frac{1}{2}}N
$$
\end{lemma}

\begin{proof} By Lemma \ref{lem:multiplier_bound1} it suffices to consider the
region $A_0$ and to show that
$$
\lb\tau+\xi^2\rb^{-\frac{1}{2}}M_0 \leq 8 \tilde{M}_0+4\tilde{N}
$$
\begin{enumerate}
\item\label{case:21} $|\xi|> 2 |\xi_1|$ and $|\xi|> 2 |\xi_2|$: In this case $|\xi_3|\leq 2|\xi|$.
In $A_0$ we have $16 \lb\tau+\xi^2\rb \geq \lb\xi\rb^2$,
since $\lb\tau+\xi^2\rb \geq
\lb\tau_1+\xi_1^2\rb,\lb\tau_2+\xi_2^2\rb,\lb\tau_3-\xi_3^2\rb$
which implies $$ 8 \lb\tau+\xi^2\rb^{\frac{1}{2}} \geq
\lb\tau_1+\xi_1^2\rb^{\delta}
\lb\tau_2+\xi_2^2\rb^{\delta}\lb\tau_3-\xi_3^2\rb^{\delta}\lb\xi\rb^{\frac{1}{2}-3\delta}
\lb\xi_3\rb^{\frac{1}{2}-3\delta}$$
\item\label{case:22} $|\xi|\leq 2 |\xi_1|$ and $|\xi|\leq 2 |\xi_2|$:
In this case we have $|\xi_3|\leq 4 \max\{|\xi_1|,|\xi_2|\}$ and $|\xi|\leq 2
\min\{|\xi_1|,|\xi_2|\}$, which shows $|\tilde{M}|\leq 4 \tilde{N}$.
\item\label{case:23} $|\xi|> 2 |\xi_1|$ and $|\xi|\leq 2 |\xi_2|$: Here
  $\xi\not=0$ and without loss we may assume $\xi_3 \not=0$, since otherwise $\tilde{M}=0$.
We have $$|\xi| \leq 2|\xi-\xi_1| \text{ and }
|\xi||\xi-\xi_2|\leq 2 \lb(\xi-\xi_1)(\xi-\xi_2)\rb$$
In the subregion where $|\xi_1|\leq |\xi-\xi_2|$ we have
$$|\xi_3|\leq
|\xi-\xi_2|+|\xi_1|\leq 2|\xi-\xi_2|$$
and therefore
$$\lb\xi\rb\lb\xi_3\rb\leq 2|\xi||\xi_3|
\leq 4|\xi||\xi-\xi_2|\leq 8 \lb(\xi-\xi_1)(\xi-\xi_2)\rb
$$
which is bounded by $32\lb\tau+\xi^2\rb$, since we are in region $A_0$.
Then,
$$8\lb\tau+\xi^2\rb^{\frac{1}{2}} \geq
\lb\xi\rb^{\frac{1}{2}-3\delta}\lb\xi_3\rb^{\frac{1}{2}-3\delta}
\lb\tau_1+\xi_1^2\rb^{\delta}\lb\tau_2+\xi_2^2\rb^{\delta}\lb\tau_3-\xi_3^2\rb^{\delta}$$
which proves $$\lb\tau+\xi^2\rb^{-\frac{1}{2}}M_0 \leq 8 \tilde{M}_0$$
In the subregion where $|\xi_1|> |\xi-\xi_2|$ we have
$|\xi_3|\leq 2|\xi_1|$
and we arrive at $\tilde{M} \leq 4 \tilde{N}$.
\item\label{case:24} $|\xi|\leq 2 |\xi_1|$ and $|\xi|> 2 |\xi_2|$:
By the symmetry of $\tilde{M}$ in $\xi_1,\xi_2$ we find the same estimate as in case \ref{case:3}.
\end{enumerate}
\end{proof}

\begin{theorem}\label{thm:trilinear_est2}
There exists $c,\eps>0$, such that for $T \in (0,1]$ and $u_j \in \test$
with $\supp(u_j) \subset \{(t,x) \mid
|t|\leq T\}$, $j=1,2,3$, we have 
\begin{equation}
  \label{eq:trilinear_est2}
  \left\|u_1 u_2 \partial_x u_3\right\|_{Y_
{\frac{1}{2},-1}}
\leq c T^{\eps} \|u_1\|_{X_{\frac{1}{2},\frac{1}{2}}}
\|u_2\|_{X_{\frac{1}{2},\frac{1}{2}}}\|u_3\|_{X^-_{\frac{1}{2},\frac{1}{2}}}
\end{equation}
\end{theorem}

\begin{proof}
We use the notation from the proof of Theorem \ref{thm:trilinear_est1}.
With the Fourier multiplier $\tilde{M}$ defined in Lemma
\ref{lem:multiplier_bound1}
we rewrite the left hand side as
\begin{equation*}
\left\|u_1 u_2 \partial_x u_3\right\|_{Y_{\frac{1}{2},-1}}=\left\|
\int_\ast
\sum_\ast \tilde{M}(\tau_1,\tau_2,\tau_3,\xi_1,\xi_2,\xi_3)
\prod_{j=1}^3 f_j(\tau_j,\xi_j)
\right\|_{l^2_{\xi}L^1_{\tau}}
\end{equation*}
By the estimate \eqref{eq:multiplier_bound2}
we successively replace $\tilde{M}$ by
$\tilde{M}_0,\tilde{M}_1,\tilde{M}_2,\tilde{M}_3$ and $\tilde{N}$.\\
\noindent\emph{Estimate for $\tilde{M}_0$:}
We observe that by the Cauchy-Schwarz inequality we have for fixed $\xi$
\begin{equation}\label{eq:cs}
\left\|\lb\tau+\xi^2\rb^{-\frac{1}{2}-\delta'}\phi(\tau,\xi)\right\|_{L^1_{\tau}}
\leq \left(\int\lb\tau\rb^{-1-2\delta'}\,
  d\tau\right)^{\frac{1}{2}}\left\| \phi(\cdot,\xi) \right\|_{L^2_{\tau}}
\end{equation}
for $\delta'>0$. Now, for fixed $\xi_1,\xi_2,\xi_3$ and $\xi=\xi_1+\xi_2+\xi_3$
\begin{align*}
&\left\|
\int_\ast
\tilde{M}_0(\tau_1,\tau_2,\tau_3,\xi_1,\xi_2,\xi_3)
\prod_{j=1}^3 f_j(\tau_j,\xi_j)
\right\|_{L^1_{\tau}}\\
=&\lb\xi\rb^{-\frac{1}{2}+3\delta}\left\|
\int_{\tau=\tau_1+\tau_2+\tau_3}
\prod_{j=1}^2
\frac{f_j(\tau_j,\xi_j)}{\lb\xi_j\rb^{\frac{1}{2}}\lb\tau_j+\xi_j^2\rb^{\frac{1}{2}+\delta}}
\;\frac{f_3(\tau_3,\xi_3)}{\lb\xi_3\rb^{\frac{1}{2}-3\delta}\lb\tau_3-\xi_3^2\rb^{\frac{1}{2}+\delta}}
\right\|_{L^1_{\tau}}\\
\leq& c \lb\xi\rb^{-\frac{1}{2}+3\delta} \prod_{j=1}^2
\left\|\frac{f_j(\tau_j,\xi_j)}{\lb\xi_j\rb^{\frac{1}{2}}
\lb\tau_j+\xi_j^2\rb^{\delta/2}}\right\|_{L^2_{\tau}}
\;\left\|\frac{f_3(\tau_3,\xi_3)}{\lb\xi_3\rb^{\frac{1}{2}-3\delta}\lb\tau_3-\xi_3^2\rb^{\delta/2}}
\right\|_{L^2_{\tau}}
\end{align*}
by Young's inequality and \eqref{eq:cs} with $\delta'=\delta/2$. With
$g_j(\tau_j,\xi_j)=f_j(\tau_j,\xi_j)\lb\tau_j+\xi_j^2\rb^{-\delta/2}$ and
$g_3(\tau_3,\xi_3)=f_3(\tau_3,\xi_3)\lb\tau_3-\xi_3^2\rb^{-\delta/2}$ we have
\begin{align*}
&\left\|\int_\ast
\sum_\ast \tilde{M}_0(\tau_1,\tau_2,\tau_3,\xi_1,\xi_2,\xi_3)
\prod_{j=1}^3 g_j(\tau_j,\xi_j)
\right\|_{l^2_{\xi}L^1_{\tau}}\\
\leq& c
\left\|\lb\xi\rb^{-\frac{1}{2}+3\delta}\sum_{\xi=\xi_1+\xi_2+\xi_3}
\lb\xi_1\rb^{-\frac{1}{2}}
\lb\xi_2\rb^{-\frac{1}{2}}
\lb\xi_3\rb^{-\frac{1}{2}+3\delta}
\prod_{j=1}^3\left\|g_j(\cdot,\xi_j)\right\|_{L^2_{\tau}}
\right\|_{l^2_{\xi}}
\end{align*}
An application of H\"older's and Young's inequalities, choosing $\delta=1/24$, gives the
upper bound
\begin{align*}
&c
\left\|\sum_{\xi=\xi_1+\xi_2+\xi_3}
\lb\xi_1\rb^{-\frac{1}{2}}
\lb\xi_2\rb^{-\frac{1}{2}}
\lb\xi_2\rb^{-\frac{3}{8}}
\prod_{j=1}^3\left\|g_j(\cdot,\xi_j)\right\|_{L^2_{\tau}}
\right\|_{l^4_{\xi}}\\
\leq &c \prod_{j=1}^3\left\|\lb\xi_j\rb^{-\frac{3}{8}}
\left\|g_j(\cdot,\xi_j)\right\|_{L^2_{\tau}}\right\|_{l^{4/3}_{\xi}}
\leq c\prod_{j=1}^3\left\|g_j\right\|_{L^2_{\tau}l^{2}_{\xi}}\\
\leq & c 
\|u_1\|_{X_{\frac{1}{2},\frac{23}{48}}}
\|u_2\|_{X_{\frac{1}{2},\frac{23}{48}}}\|u_3\|_{X^-_{\frac{1}{2},\frac{23}{48}}}
\end{align*}
which finally proves that
\begin{align*}
\left\|
\int_\ast
\sum_\ast \tilde{M}_0(\tau_1,\ldots,\xi_3)
\prod_{j=1}^3 f_j(\tau_j,\xi_j)
\right\|_{l^2_{\xi}L^1_{\tau}}
 \leq c T^{\eps} \|u_1\|_{X_{\frac{1}{2},\frac{1}{2}}}
\|u_2\|_{X_{\frac{1}{2},\frac{1}{2}}}\|u_3\|_{X^-_{\frac{1}{2},\frac{1}{2}}}
\end{align*}
\noindent\emph{Estimate for $\tilde{M}_1,\tilde{M}_2,\tilde{M}_3$ and
  $\tilde{N}$:}
We show that the estimates from the proof of Theorem \ref{thm:trilinear_est1} are strong enough to treat
these terms, too.
Indeed, an application of \eqref{eq:cs} implies
\begin{align*}
&\left\|
\int_\ast
\sum_\ast \tilde{M}_1(\tau_1,\tau_2,\tau_3,\xi_1,\xi_2,\xi_3)
\prod_{j=1}^3 f_j(\tau_j,\xi_j)
\right\|_{l^2_{\xi}L^1_{\tau}}\\
&\leq c
\left\|\lb\tau+\xi^2\rb^{\frac{1}{8}}\int_\ast
\sum_\ast M_1(\tau_1,\tau_2,\tau_3,\xi_1,\xi_2,\xi_3)
\prod_{j=1}^3 f_j(\tau_j,\xi_j)
\right\|_{l^2_{\xi}L^2_{\tau}}= c m_1
\end{align*}
where $m_1$ is defined in \eqref{eq:m_1} and is bounded according to
\eqref{eq:bound_m_1}.
The same reasoning applies to $\tilde{M}_2,\tilde{M}_3$ and $\tilde{N}$, where
we use the bounds established in \eqref{eq:bound_m_2}, \eqref{eq:bound_m_3} and
\eqref{eq:bound_n}.
\end{proof}
The next Lemma contains an auxiliary estimate, which will be used for
polynomial terms in the nonlinearity.
This suffices for our purposes, but it is far from optimal, see \cite{Bo93}.
\begin{lemma}\label{lem:pol_est}
For $\delta>0$ there exists $c,\eps>0$, such that for $T \in (0,1]$ and $u_j \in
\test $ with $\supp(u_j) \subset \{(t,x) \mid
|t|\leq T\}$, $j=1,\ldots, 5$, we have
\begin{equation}\label{eq:quint_pol_est}
\Big\|\prod_{j=1}^5 u_j\Big\|_{X_{\frac{1}{2},-\frac{3}{8}-\delta}}\leq cT^\eps
\|u_1\|_{X^-_{\frac{1}{2},\frac{1}{2}}}\|u_2\|_{X^-_{\frac{1}{2},\frac{1}{2}}}
\prod_{j=3}^5\|u_j\|_{X_{\frac{1}{2},\frac{1}{2}}}
\end{equation}
and
\begin{equation}\label{eq:tri_pol_est}
\Big\|\prod_{j=1}^3u_j\Big\|_{X_{\frac{1}{2},-\frac{3}{8}-\delta}}
\leq cT^\eps
\|u_1\|_{X^-_{\frac{1}{2},\frac{1}{2}}}\|u_2\|_{X_{\frac{1}{2},\frac{1}{2}}}
\|u_3\|_{X_{\frac{1}{2},\frac{1}{2}}}
\end{equation}
\end{lemma}
\begin{proof}
As in the previous proofs it suffices to consider $\mathcal{F}u_j \geq 0$.
For $\xi=\sum\limits_{k=1}^5 \xi_k$ we have
$
\lb\xi\rb^{\frac{1}{2}}\leq c \sum\limits_{k=1}^5 \lb\xi_k\rb^\frac{1}{2}
$
which implies
$$
\Big\|\prod_{j=1}^5 u_j\Big\|_{X_{\frac{1}{2},-\frac{3}{8}-\delta}}\leq c
\sum_{k=1}^5\Big\|J^{\frac{1}{2}}u_k\prod_{\genfrac{}{}{0pt}{}{j=1}{j\not=k}}^5 u_j\Big\|_{X_{0,-\frac{3}{8}-\delta}}
$$
Each of the five terms can be estimated, using the dual Strichartz estimate
\eqref{eq:strichartz_dual} as follows
\begin{align*}
\Big\|J^{\frac{1}{2}}u_k\prod_{\genfrac{}{}{0pt}{}{j=1}{j\not=k}}^5
u_j\Big\|_{X_{0,-\frac{3}{8}-\delta}}&\leq c\Big\|J^{\frac{1}{2}}u_k\prod_{\genfrac{}{}{0pt}{}{j=1}{j\not=k}}^5
u_j\Big\|_{L^{\frac{4}{3}}_{t}L^{\frac{4}{3}}_{x}}\\
& \leq c \|J^{\frac{1}{2}}u_k\|_{L^2_{t}L^2_{x}}\prod_{\genfrac{}{}{0pt}{}{j=1}{j\not=k}}^5
\|u_j\|_{L^{16}_{t}L^{16}_{x}}\\
&\leq cT^\eps\|u_1\|_{X^-_{\frac{1}{2},\frac{1}{2}}}\|u_2\|_{X^-_{\frac{1}{2},\frac{1}{2}}}
\prod_{j=3}^5\|u_j\|_{X_{\frac{1}{2},\frac{1}{2}}}
\end{align*}
where in the last step we used the Sobolev embedding in space and time.
The second claim follows in the same way, using the $L^8_tL^8_x$ norm instead of
the $L^{16}_tL^{16}_x$ norm on the
factors without derivatives.
\end{proof}
We put these estimates in a slightly more general form.
\begin{corollary}\label{cor:multi}
Let $s\geq \frac{1}{2}$ and $\delta>0$. There exists $c,\eps>0$, such that for $T \in (0,1]$
and $u_j \in \test $ with $\supp(u_j) \subset \{(t,x) \mid
|t|\leq T\}$, $j=1,\ldots, 5$, we have
\begin{equation}
  \label{eq:gen_tri}
  \left\|u_1 u_2 \partial_x \overline{u}_3\right\|_{Y_{s,-1} \cap X_{s,-\frac{1}{2}}}
\leq c
T^{\eps}\sum\limits_{k=1}^3\|u_k\|_{X_{s,\frac{1}{2}}}
\prod_{\genfrac{}{}{0pt}{}{j=1}{j\not=k}}^3\|u_j\|_{X_{\frac{1}{2},\frac{1}{2}}}
\end{equation}
\begin{equation}\label{eq:gen_quint_pol_est}
\Big\|\overline{u}_1\overline{u}_2\prod_{j=3}^5 u_j\Big\|_{X_{s,-\frac{3}{8}-\delta}}\leq cT^\eps
\sum\limits_{k=1}^5\|u_k\|_{X_{s,\frac{1}{2}}}
\prod_{\genfrac{}{}{0pt}{}{j=1}{j\not=k}}^5\|u_j\|_{X_{\frac{1}{2},\frac{1}{2}}}
\end{equation}
\begin{equation}\label{eq:gen_tri_pol_est}
\begin{split}
&|\mu(u_1)-\mu(u_2)|\Big\|\overline{u}_3u_4u_5\Big\|_{X_{s,-\frac{3}{8}-\delta}}\\
\leq
&cT^\eps
\|u_1-u_2\|_{Z_0}(\|u_1\|_{Z_0}+\|u_2\|_{Z_0})
\sum\limits_{k=3}^5\|u_k\|_{X_{s,\frac{1}{2}}}
\prod_{\genfrac{}{}{0pt}{}{j=3}{j\not=k}}^5\|u_j\|_{X_{\frac{1}{2},\frac{1}{2}}}
\end{split}
\end{equation}
and
\begin{equation}\label{eq:gen_psi_est}
\begin{split}
&\Big\|(\psi(u_1)-\psi(u_2))u_3\Big\|_{X_{s,0}}\\
\leq&  cT^\eps
(1+\|u_1\|_{X_{\frac{1}{2},\frac{1}{2}}\cap Z_0}+\|u_2\|_{X_{\frac{1}{2},\frac{1}{2}}\cap Z_0})^3
\|u_1-u_2\|_{X_{\frac{1}{2},\frac{1}{2}}\cap Z_{0}}\|u_3\|_{X_{s,\frac{1}{2}}}
\end{split}
\end{equation}
\end{corollary}
\begin{proof}
We observe that $\|\overline{u}\|_{X_{s,b}}=\|u\|_{X^-_{s,b}}$ and
$$
\lb\xi\rb^{s}\leq c \sum\limits_{k=1}^l \lb\xi_k\rb^s
\quad,\text{ for }\xi=\sum\limits_{k=1}^l \xi_k \text{ and } s\geq 0
$$
Furthermore, by the embedding $Z_0 \hookrightarrow C(\R,L^2(\T))$
$$|\mu(u_1)-\mu(u_2)|\leq c \|u_1-u_2\|_{Z_0}(\|u_1\|_{Z_0}+\|u_2\|_{Z_0})$$
and by \eqref{eq:psi_est}
\begin{equation*}
\|\psi(u)-\psi(v)\|_{L^4_T}\leq c T^{\eps}
(1+\|u_1\|_{X_{\frac{1}{2},\frac{1}{2}}\cap Z_0}+\|u_2\|_{X_{\frac{1}{2},\frac{1}{2}}\cap Z_0})^3
\|u_1-u_2\|_{X_{\frac{1}{2},\frac{1}{2}}\cap Z_{0}}
\end{equation*}
Using this, the corollary
follows from \eqref{eq:trilinear_est1},
\eqref{eq:trilinear_est2}, \eqref{eq:quint_pol_est} and \eqref{eq:tri_pol_est}.
\end{proof}

\section{The gauge equivalent Cauchy problem}\label{sect:thm_gauge_eq_problem}

\begin{theorem}\label{thm:main2}
Let $s\geq \frac{1}{2}$. There exists a non-increasing function $T:
(0,\infty)\to (0,\infty)$, such that for $v_0 \in H^s(\T)$ and
$T=T(\|v_0\|_{H^{\frac{1}{2}}(\T)})$ there exists a solution
$$v \in Z^{T}_{s} \subset C([-T,T],H^s(\T))
$$
of the Cauchy problem
\begin{equation}\label{eq:gauge_eq_cp}
\begin{split}
\partial_t v-i\partial^2_{x}v
&=-v^2\partial_x\overline{v}+\tfrac{i}{2}|v|^4v-i\mu(v)|v|^2v+i\psi(v)v \quad\text{in } (-T,T)\times \T\\
v(0)&=v_0
\end{split}
\end{equation}
where
$\mu(v)=\frac{1}{2\pi}\|v(0)\|^2_{L^2(\T)}$ and
$$
\psi(v)(t)=\frac{1}{2\pi}\int_0^{2\pi}
2\Imag(\overline{v}_{x}v)(t,\theta)-\frac{1}{2}|v|^4(t,\theta)\,d\theta
+\mu(v)^2
$$
This solution is unique in $Z^{T}_{\frac{1}{2}}$.
Moreover, for any $r>0$ there exists $T=T(r)$, such that with
$$B_{r}=\{v_0 \in H^{s}(\T)\mid \|v_0\|_{H^\frac{1}{2}(\T)}< r\}$$
the flow map 
$$
\tilde{F}: H^s(\T) \supset B_{r}
\to C\big([-T,T],H^{s}(\T)\big) \quad,\; v_0 \mapsto v
$$
is Lipschitz continuous.
\end{theorem}
\begin{remark}\label{rem:ill}
We remark that Theorem \ref{thm:main2} extends to nonlinear
terms of the type $\overline{u}^{k}\partial_x \overline{u}$ by Gr\"unrock's result
\cite{G00}. On the other hand, Christ
\cite{C03} proved a strong ill-posedness result for the nonlinearities
$u^k\partial_x u$, for every $k \in \N$.

\end{remark}
As in the case of the real line, we show that below $s=\frac{1}{2}$ it is not possible to prove similar
estimates on the tri-linear term which contains the derivate.
\begin{theorem}\label{thm:ill_posedness_conj_dnls}
Let $s<\frac{1}{2}$ and $T>0$. There does not exist a normed space
$Z_T \hookrightarrow C([-T,T],H^s(\T))$, such that
$$
\|W(t)u_0\|_{Z_T} \leq c\|u_0\|_{H^s(\T)}
$$
and
$$
\left\|\int_0^t W(t-t')\left(u^2 \partial_x \overline{u}\right)(t')\, dt'\right\|_{Z_T}
\leq c
\|u\|^3_{Z_T}
$$
hold.
\end{theorem}
The proof of local well-posedness of the gauge equivalent
problem will be a straightforward application of the contraction mapping
principle, cp. \cite{Bo93,CKSTT03}.
We define for $v \in \test$ 
$$
N(v)=-v^2\partial_x\overline{v}+\tfrac{i}{2}|v|^4v-i\mu(v)|v|^2v+i\psi(v)v
$$
where $\mu(v)(t)=\frac{1}{2\pi}\|v(t)\|^2_{L^2}$ and
$$\psi(v)(t)=\frac{1}{2\pi}\int_0^{2\pi}
2\Imag(\overline{v}_{x}v)(t,\theta)-\frac{1}{2}|v|^4(t,\theta)\,d\theta
+\frac{1}{4\pi^2}\|v(t)\|^4_{L^2}
$$
and $N_T(v)=N(\chi_T v)$ as well as
\begin{equation}\label{eq:nonl_op}
\Phi_{T}(v)(t)=\chi(t)\int_0^t W(t-t') N_T(v)(t')\, dt'
\end{equation}
We recall the definition of the space
$$Z_s=X_{s,\frac{1}{2}} \cap Y_{s,0}$$
see \eqref{eq:X_s,b}, \eqref{eq:Y_s} and \eqref{eq:Z}.
By Corollary \ref{cor:multi}, the embedding \eqref{eq:x_in_y} and the linear estimate
\eqref{eq:lin_inhom}, we may extend $\Phi_T$ uniquely to
$$\Phi_{T}: Z_s\to Z_s$$
for all $s\geq \frac{1}{2}$. We also have $$\Phi_T \big|_{[-T,T]}: \, Z_s^T \to
Z_s^T$$
since it only depends on $v\big|_{[-T,T]}$.

\begin{proof}[Local Existence]
Our aim is to find a solution $v \in Z_s$ of
\begin{equation*}
v=\chi W(\cdot)v_0+\Phi_{T}(v)
\end{equation*}
For $v_0 \in H^s(\T)$
we use again the estimates from Corollary \ref{cor:multi} and
\eqref{eq:x_in_y}, \eqref{eq:lin_hom} and \eqref{eq:lin_inhom} as well as Lemma
\ref{lem:cutoffs} to show that
there exists $c,\eps>0$, such that
\begin{equation*}
\|\chi W(\cdot)v_0+\Phi_{T}(v)\|_{Z_s}\leq c\|v_0\|_{H^s}+cT^\eps (1+\|v\|_{Z_s})^3\|v\|^2_{Z_s}
\end{equation*}
and
\begin{equation*}
\|\Phi_{T}(v_1)-\Phi_{T}(v_2)\|_{Z_s}\leq cT^\eps
(1+\|v_1\|_{Z_s}+\|v_2\|_{Z_s})^3(\|v_1\|_{Z_s}+\|v_2\|_{Z_s})\|v_1-v_2\|_{Z_s}
\end{equation*}
Then, for all $v_0 \in H^s$ with $\|v_0\|_{H^s} \leq r$ and $R=2cr$ and $T>0$ so small
that $T \leq (4c^2r(1+4cr)^3)^{-\frac{1}{\eps}}$ we see that
\begin{equation*}
v \mapsto \chi W(\cdot)v_0+\Phi_{T}(v)
\end{equation*}
maps the closed ball $B_R \subset Z_s$ to itself and is a strict contraction.
This shows the existence of a solution $v\in B_R \subset Z_s$. By restriction to
the interval $[-T,T]$ we found a solution $v \in Z^T_s\subset C([-T,T],H^s(\T))$ of
\begin{equation}\label{eq:gauge_int_eq}
v(t)=W(t)v_0+\Phi_{T}(v)(t) \quad, \; t \in [-T,T]
\end{equation}
\end{proof}

\begin{proof}[Uniqueness]
Assume that $v_1,v_2 \in Z^T_{\frac{1}{2}}$ are two solutions of
\eqref{eq:gauge_int_eq}, such that
$$
T':=\sup\{t \in [0,T] \mid v_1(t)=v_2(t)\}<T
$$
and we define $w_j(t)=\tilde{v}_j(T'+t)$, $j=1,2$ for extensions $\tilde{v}_j$
of $v_j$. By approximation we see
$$
w_1(t)-w_2(t)=\Phi_T(w_1)(t)-\Phi_T(w_2)(t) \quad -T' \leq t \leq T-T'
$$
Choosing $\delta>0$ small enough, we arrive at
$$
\|\chi_\delta (w_1-w_2)\|_{Z_\frac{1}{2}} \leq c
\delta^{\eps}(1+\|w_1\|_{Z_\frac{1}{2}}
+\|w_2\|_{Z_\frac{1}{2}})^4\|\chi_\delta (w_1-w_2)\|_{Z_\frac{1}{2}}
$$
which forces $w_1(t)=w_2(t)$ for $|t| \leq \delta$ and therefore contradicts the
definition of $T'$. The same argument applies in the interval $[-T,0]$.
\end{proof}

\begin{proof}[Local Lipschitz continuity of the flow]
Let $v_0,w_0 \in H^s(\T)$ with $\|v_0\|_{H^s},\|w_0\|_{H^s} \leq r$. Let $v,w
\in Z^T_s$ be two solutions
of \eqref{eq:gauge_int_eq} with $v(0)=v_0$ and $w(0)=w_0$ with
extensions $\tilde{v},\tilde{w}$ constructed in part 1 of the proof. Then,
$\|\tilde{v}\|_{Z_s},\|\tilde{w}\|_{Z_s} \leq 2cr$ and
\begin{equation*}
\|\tilde{v}-\tilde{w}\|_{Z_s}\leq c\|v_0-w_0\|_{H^s}+cT^\eps
4c^2r(1+4cr)^3\|\tilde{v}-\tilde{w}\|_{Z_s}
\end{equation*}
and the choice of $T$ from part 1 guarantees
$$
\|\tilde{v}-\tilde{w}\|_{Z_s}\leq 2c\|v_0-w_0\|_{H^s}
$$ 
and by restriction
$$
\|v-w\|_{C([-T,T],H^s)}\leq 2c\|v_0-w_0\|_{H^s}
$$
\end{proof}

\begin{proof}[Time of existence]
Finally, the standard iteration argument, using the estimates from Corollary \ref{cor:multi},
shows that the maximal time of existence $T>0$ depends only on
$\|v_0\|_{H^\frac{1}{2}}$.
\end{proof}

Finally, we remark that the counterexamples from \cite{Tak99} also show the
optimality of our tri-linear estimate:\\
\begin{proof}[Proof of Theorem \ref{thm:ill_posedness_conj_dnls}]
We follow the general idea from \cite{MST01}.
Let $n \in \N$ and $u^{(n)}_0:=n^{-s}e^{inx}$. Then, $\|u^{(n)}_0\|_{H^s}=c$ and
$$
\int_0^t W(-t')\left((W(t')u^{(n)}_0)^2\partial_x \overline{W(t')u^{(n)}_0}\right) \,
dt'=-it n^{-3s}ne^{inx}
$$
which shows that
$$
\left\|\int_0^t W(t-t')\left((W(t')u^{(n)}_0)^2\partial_x
    \overline{W(t')u^{(n)}_0}\right) \, dt'\right\|_{H^s} \geq c|t|n^{1-2s}
$$
If the linear and tri-linear estimates in a space $Z_T \hookrightarrow
C([-T,T],H^s(\T))$ were true, there would exist $c>0$ such that
$
|t|n^{1-2s}\leq c
$ for all $n \in \N$,  which is a contradiction for $s<\frac{1}{2}$.
\end{proof}

\section{Proof of Theorem \ref{thm:main} and Corollary \ref{cor:global}}\label{sect:proof_thm_main}
In this section, we will use the solutions of \eqref{eq:gauge_eq_cp} constructed in the previous section 
to prove Theorem \ref{thm:main}, similar to \cite{Tak99,G05}.

\begin{proof}[Existence]
We fix $s\geq \frac{1}{2}$ and let $u_0 \in
H^s(\T)$ with $\mu:=\frac{1}{2\pi}\|u_0\|^2_{L^2}$. Then, we define
$v_0:=\mathcal{G}(u_0) \in H^s(\T)$, see Lemma
\ref{lem:gauge_trafo_est}. According to Theorem \ref{thm:main2}, there exists a
unique solution $v \in Z^{T}_s \subset C([-T,T],H^s(\T))$ of
\eqref{eq:gauge_int_eq}. Now, we claim
that $u:=\mathcal{G}^{-1}(v) \in \mathcal{X}_s^T \subset C([-T,T],H^s(\T))$ solves
\begin{equation}\label{eq:solution_DNLS}
u(t)=W(t)u_0+\int_{0}^t W(t-t') \partial_x(|u|^2u)(t')\, dt', \quad t \in (-T,T)
\end{equation}
For smooth functions this follows from Lemma \ref{lem:gauge_eq_prob}.
Let $u^{(n)}_0 \in C^{\infty}$ with $u^{(n)}_0 \to u_0$ in $H^s$ and
$\|u^{(n)}_0\|_{L^2}=\|u_0\|_{L^2}$. Moreover, let $v^{(n)}\in Z^{T}_s $ be the solution of
\eqref{eq:gauge_int_eq} with initial data $\mathcal{G}(u^{(n)}_0)$ and
$u^{(n)}:=\mathcal{G}^{-1}(v^{(n)})$.
Then,
\begin{align*}
& \sup_{t \in (-T,T)}\left\|\int_{0}^t W(t-t')
  \partial_x(|u|^2u-|u^{(n)}|^2u^{(n)})(t')\, dt'\right\|_{H^{-1}}\\
\leq & c
  (\|u\|^2_{L^{\infty}_TH^{\frac{1}{2}}}+\|u^{(n)}\|^2_{L^{\infty}_TH^{\frac{1}{2}}})
\|u-u^{(n)}\|_{L^1_TL^2_x}
\end{align*}
Because $\mathcal{G}$ is continuous in $H^s$, $\mathcal{G}(u^{(n)}_0)\to v_0$ and due to
the continuity of the flow map of \eqref{eq:gauge_int_eq} we have
$v^{(n)}\to v$ in $C([-T,T],H^s)$. Since also $\mathcal{G}^{-1}$ is continuous, the above term tends to zero. This
shows that $u$ solves
\eqref{eq:solution_DNLS}
because obviously also the linear part converges in $H^s(\T)$.
\end{proof}

\begin{proof}[Uniqueness]
Let $u_1,u_2 \in \mathcal{X}_s^T$ be two solutions of
\eqref{eq:solution_DNLS} with $u_1(0)=u_2(0)$,
such that $\mathcal{G} (u_j) \in Z^T_{\frac{1}{2}}$ solve \eqref{eq:gauge_int_eq} with
the same initial datum. By the uniqueness of the
solutions to \eqref{eq:gauge_int_eq} we have $\mathcal{G}
(u_1)=\mathcal{G} (u_2)$ and therefore $u_1=u_2$.

We now prove that the hypothesis that $\mathcal{G} (u_j)$ solve \eqref{eq:gauge_int_eq} is
fulfilled if $u_j$ are limits of smooth solutions in $\mathcal{X}_{\frac{1}{2}}^T$, say $u_j^{(n)}
\in C([-T,T],H^3(\T)) \cap C^1([-T,T],H^1(\T))$ such that
$\|\mathcal{G}(u_j^{(n)})-\mathcal{G}(u_j)\|_{Z^T_{\frac{1}{2}}} \to 0$. By Lemma \ref{lem:gauge_eq_prob}
$\mathcal{G}(u^{(n)}_j) \in Z^T_s$ solve \eqref{eq:gauge_int_eq}. Moreover,
$\mathcal{G}(u_j^{(n)}(0))\to \mathcal{G}(u_j(0)) \in H^\frac{1}{2}$. There exists a unique solution $v
\in Z^T_{\frac{1}{2}}$ to \eqref{eq:gauge_int_eq} with
$v(0)=\mathcal{G}(u_j(0))$ and due to the continuity of the flow $\tilde{F}$ it
follows $\mathcal{G}(u^{(n)}_j) \to v $, which implies that $\mathcal{G}
(u_j)=v$ is a solution to \eqref{eq:gauge_int_eq}.
\end{proof}

\begin{proof}[Continuity of the flow]
Since the flow map to \eqref{eq:solution_DNLS} $F:H^s(\T) \to C([-T,T],H^s(\T))$ results from conjugating the flow
map to \eqref{eq:gauge_int_eq} $\tilde{F}:H^s(\T) \to C([-T,T],H^s(\T))$ with the gauge transformation
$\mathcal{G}$, i.e. $F=\mathcal{G}^{-1} \circ \tilde{F} \circ \mathcal{G}$, its continuity
properties follow from the local Lipschitz continuity of $\tilde{F}$ and Lemma
\ref{lem:gauge_trafo_est}.
\end{proof}

\begin{proof}[Global existence]
It suffices to prove an a priori bound for smooth solutions.
By Lemma \ref{lem:energy_cons} and the Sobolev embedding we have
\begin{equation}\label{eq:c}
\|\partial_x u(t)\|_{L^2(\T)}^2+\frac{3}{2}\Imag \int_0^{2\pi} |u|^2u
\partial_x\overline{u}(t) \, dx +\frac{1}{2}\|u(t)\|^6_{L^6(\T)}
\leq c(1+ \|u_0\|_{H^1(\T)})^6
\end{equation}
Now, we use the Gagliardo-Nirenberg inequality
\begin{equation}\label{eq:GN}
\|u(t)\|^3_{L^6(\T)} \leq \|u(t)\|^2_{L^2(\T)}\left( \|\partial_x u(t)\|_{L^2(\T)}+\frac{1}{2\pi}\|u(t)\|_{L^2(\T)}\right)
\end{equation}
see Appendix \ref{app_sect:GN}, and estimate
\begin{equation*}
\frac{3}{2}\Imag \int_0^{2\pi} |u|^2u
\partial_x\overline{u} (t)\, dx
\geq -\frac{3}{2}\|u(t)\|^2_{L^2}\left(\|\partial_x
  u(t)\|_{L^2}+\frac{1}{2\pi}\|u(t)\|_{L^2}\right)\|\partial_x u(t)\|_{L^2}
\end{equation*}
Then, using this in \eqref{eq:c} we have for $\|u(t)\|_{L^2}\leq \delta$
$$
(1-\frac{3}{2}\delta^2)\|\partial_x u(t)\|_{L^2}^2-\frac{3}{4\pi}\delta^3\|\partial_x u(t)\|_{L^2}
\leq c(1+ \|u_0\|_{H^1(\T)})^6
$$
which shows for $\delta< \sqrt{\frac{2}{3}}$ that there exists $c(\delta)>0$
such that
$$
\|\partial_x u(t)\|_{L^2}^2
\leq c(\delta)(1+ \|u_0\|_{H^1(\T)})^6
$$
This estimate, together with the $L^2$ conservation law from Lemma
\ref{lem:l2_cons} shows that for
$\|u_0\|_{L^2}\leq \delta$ there exists $C(\delta)>0$ such that
$$
\|u(t)\|_{H^1(\T)}
\leq C(\delta)(1+ \|u_0\|_{H^1(\T)})^3
$$
\end{proof}

\begin{remark}\label{rem:global}
The proof shows that it suffices to choose $\delta<\sqrt{\frac{2}{3|\lambda|}}$.
By following the idea of Hayashi and Ozawa \cite{HO94}, using the gauge transform
together with sharp versions of the Gagliardo-Nirenberg inequality
we expect that this can be improved, but our aim here is to give a short proof of the
qualitative result.
\end{remark}

\begin{appendix}
\section{Proof of estimate \eqref{eq:exp_est}}\label{app_sect:exp_est}
We prove that for all $s\geq 0$ there exists $c>0$, such that for $f,g,h \in
H^s(\T)$ we have
$$
\left\|(e^{\pm i \mathcal{I}(f)}-e^{\pm i
  \mathcal{I}(g)})h\right\|_{H^s}\leq c
e^{c\|f\|^2_{H^s}+c\|g\|^2_{H^s}}(\|f\|_{H^s}+\|g\|_{H^s})\|f-g\|_{H^s}\|h\|_{H^s}
$$
To simplify the notation
we only consider the plus sign since the same argument works with the minus
sign. Moreover, it suffices to consider smooth $f,g,h$ and we start with the
case $s>0$.
We will exploit the Sobolev multiplication law
$$\|fg\|_{H^\alpha} \leq c
\|f\|_{H^\alpha}\|g\|_{H^{\beta}}\quad,
\beta=\begin{cases}\alpha & ,\alpha>\frac{1}{2}\\
\frac{1}{2}+\eps & ,\text{otherwise}
\end{cases}
$$
We write
\begin{equation}\label{eq:series}
(e^{i \mathcal{I}(f)}-e^{i
  \mathcal{I}(g)})h=ih(\mathcal{I}(f)-\mathcal{I}(g))\sum_{k=1}^{\infty}\frac{1}{k!}\sum_{j=0}^{k-1}
(i\mathcal{I}(f))^j(i\mathcal{I}(g))^{k-1-j}
\end{equation}
Let $s'=\max\{s,\frac{1}{2}+\eps\}$ for some $0<\eps<\frac{1}{2}$ to be chosen
later. Then, the $H^s$ norm of the expression \eqref{eq:series} is bounded by
$$
\|h\|_{H^s}\|\mathcal{I}(f)-\mathcal{I}(g)\|_{H^{s'}}\sum_{k=1}^{\infty}\frac{1}{k!}\sum_{j=0}^{k-1}
(c\|\mathcal{I}(f)\|_{H^{s'}})^j(c\|\mathcal{I}(g)\|_{H^{s'}})^{k-1-j}
$$
Now, we observe that
$$
\sum_{k=1}^{\infty}\frac{1}{k!}\sum_{j=0}^{k-1}
(c\|\mathcal{I}(f)\|_{H^{s'}})^j(c\|\mathcal{I}(g)\|_{H^{s'}})^{k-1-j}
\leq e^{c \|\mathcal{I}(f)\|_{H^{s'}}+c \|\mathcal{I}(g)\|_{H^{s'}}}
$$
Moreover,
$$
\|\mathcal{I}(f)\|_{H^{s'}}\leq \||f|^2\|_{H^{s'-1}}+\|f\|^2_{L^2}
$$
In the case where $s\geq\frac{1}{2}+\eps$ it follows $\||f|^2\|_{H^{s'-1}}\leq
\||f|^2\|_{H^{s}}\leq c\|f\|^2_{H^{s}}$ and otherwise, with $p=\frac{1}{1-\eps}$
$$
\||f|^2\|_{H^{-\frac{1}{2}+\eps}} \leq c \||f|^2\|_{L^p}\leq c \|f\|^2_{L^{2p}} \leq c\|f\|^2_{H^{\frac{\eps}{2}}}
$$
by Sobolev embeddings. Now, choosing $\eps\leq 2s$ we have
$$
\|\mathcal{I}(f)\|_{H^{s'}}\leq c \|f\|^2_{H^{s}}
$$
Similarly, we get
$$
\|\mathcal{I}(f)-\mathcal{I}(g)\|_{H^{s'}}\leq c (\|f\|_{H^{s}}+\|g\|_{H^{s}})\|f-g\|_{H^{s}}
$$
and the claim follows for $s>0$.
Finally, for $s=0$
\begin{align*}
\left\|(e^{i \mathcal{I}(f)}-e^{i \mathcal{I}(g)})h\right\|_{L^2}& \leq
\left\|e^{ i \mathcal{I}(f)}-e^{ i \mathcal{I}(g)}\right\|_{L^\infty}\|h\|_{L^2}\\
& \leq
\left\|\mathcal{I}(f)-\mathcal{I}(g)\right\|_{L^\infty}\|h\|_{L^2}\\
&\leq 
2(\|f\|_{L^2}+\|g\|_{L^2})\|f-g\|_{L^2}\|h\|_{L^2}
\end{align*}

\section{Conservation laws}\label{app_sect:cons_q}
The results in this section are well-known in the case of the real line (cp. \cite{CH98}, Proposition
6.1.1, appendix of \cite{HO94}, or \cite{KN78}) and formally everything transfers to the periodic
setting. Nevertheless, we briefly repeat the main points for completeness of the
paper.
\begin{lemma}\label{lem:l2_cons}
If $$u \in C([-T,T],H^2(\T)) \cap C^1([-T,T],L^2(\T))$$
is a solution of \eqref{eq:pdnls} or \eqref{eq:gauge_eq_cp}, we have for $t\in (-T,T)$
$$
\diff{t} \|u(t)\|_{L^2(\T)}=0
$$
\end{lemma}
\begin{proof}
One easily shows that
$$
\diff{t}\|u(t)\|^2_{L^2(\T)}=2\Real \int_0^{2\pi} \overline{u} N(u)(t) \, dx
$$
for $N(u)=\partial_x (|u|^2 u)$ or
$N(u)=-u^2\partial_x\overline{u}+\tfrac{i}{2}|u|^4u-i\mu(u)|u|^2u+i\psi(u)u$,
respectively.
Partial integration yields
$$\Real \int_0^{2\pi} \overline{u} \partial_x (|u|^2 u)  \, dx=0$$
and
$$\Real \int_0^{2\pi} \overline{u} u^2 \partial_x\overline{u}  \, dx=0$$
Obviously, all the other terms also vanish and the $L^2$ conservation law follows.
\end{proof}

\begin{lemma}\label{lem:energy_cons}
If $$u \in C([-T,T],H^3(\T)) \cap C^1([-T,T],H^1(\T))$$
is a solution of \eqref{eq:pdnls} with $\lambda=1$, we have for $t\in (-T,T)$
$$
\diff{t} \left(\|u_x(t)\|_{L^2(\T)}^2+\frac{3}{2}\Imag \int_0^{2\pi} |u|^2u
\overline{u}_x (t)\, dx +\frac{1}{2}\|u(t)\|^6_{L^6(\T)}\right)=0$$
\end{lemma}
\begin{proof}
Firstly, using \eqref{eq:pdnls} we verify
  \begin{equation}
    \label{eq:dxu}
    \diff{t}\|u_x\|^2_{L^2}=2\Real \int_0^{2\pi}(|u|^2u)_{xx}\overline{u}_{x} \, dx
  \end{equation}
Secondly, we again exploit \eqref{eq:pdnls} and carry out all the differentiations
\begin{align}\label{eq:im}
&  \diff{t}\Imag \int_0^{2\pi}|u|^2u\overline{u}_x\, dx=\Imag
  \int_0^{2\pi}(|u|^2u)_t \overline{u}_x\, dx+\Imag
  \int_0^{2\pi} |u|^2u \overline{u}_{tx}\, dx\nonumber\\
&= 4 \Real \int_0^{2\pi}|u|^2u_x\overline{u}_{xx}\, dx+4\Imag \int_0^{2\pi}\overline{u}_x^2u^2|u|^2\, dx
\end{align}
Thirdly,
\begin{equation*}
  \diff{t}\|u\|^6_{L^6}=-6\Imag\int_0^{2\pi}|u|^4 \overline{u}u_{xx}\,dx
+6\Real\int_0^{2\pi}|u|^4\overline{u}(|u|^2u)_x\,dx
\end{equation*}
and
$$\Real\int_0^{2\pi}|u|^4\overline{u}(|u|^2u)_x\,dx=
\frac{3}{8}\int_0^{2\pi}(|u|^8)_x\,dx=0$$
Moreover, we integrate by parts and obtain
\begin{equation}
  \label{eq:l6}
  \diff{t}\|u\|^6_{L^6}=-12\Imag\int_0^{2\pi}|u|^2 u^2\overline{u}^2_x\,dx
\end{equation}
Now, combining \eqref{eq:dxu}, \eqref{eq:im} and \eqref{eq:l6} and integrating
by parts we get
\begin{align*}
&\diff{t}\left( \|u\|^2_{L^2}+\frac{3}{2}\Imag \int_0^{2\pi} |u|^2u
\overline{u}_x (t)\, dx +\frac{1}{2}\|u(t)\|^6_{L^6(\T)}\right)\\
=&6\Real \int_0^{2\pi}|u|^2u_x\overline{u}_{xx}\, dx-2\Real
\int_0^{2\pi}(|u|^2u)_{x}\overline{u}_{xx} \, dx\\
=&2\Real \int_0^{2\pi}|u|^2u_x\overline{u}_{xx}\, dx
-2\Real \int_0^{2\pi}u^2\overline{u}_x\overline{u}_{xx} \, dx=0
\end{align*}
\end{proof}
\section{Proof of the estimate \eqref{eq:GN}}\label{app_sect:GN}
Let $f,g$ be smooth and $2\pi$-periodic with $g(0)=0$. Then,
$g(x)=\int_0^x g'(y) \,dy$ and $g(x)=-\int^{2\pi}_x g'(y) \,dy$ such that
$$2\|g\|_{L^{\infty}(\T)} \leq \|g'\|_{L^1(\T)}$$
and by H\"older's inequality
$$\|fg\|_{L^2(\T)} \leq \frac{1}{2}\|f\|_{L^2(\T)} \|g'\|_{L^1(\T)}$$
By a translation $x\mapsto x+\xi$ we see that this holds for all $g$ with $g(\xi)=0$ for
some $\xi \in [0,2\pi]$. Now, let $u$ be smooth and $2\pi$-periodic and set
$f=u$ and $g=u^2-\frac{1}{2\pi}\int_0^{2\pi}u^2(y) \,dy$. Then,
$$\left\|u(u^2-\frac{1}{2\pi}\int_0^{2\pi}u^2(y) \, dy)\right\|_{L^2(\T)} \leq \|u\|_{L^2(\T)} \|u
u'\|_{L^1(\T)}\leq \|u\|^2_{L^2(\T)} \|u'\|_{L^2(\T)}$$
and the estimate \eqref{eq:GN} follows.
\end{appendix}
\bibliographystyle{hplain}
\bibliography{literatur}\label{sect:refs}

\end{document}